\documentclass[a4paper,11pt,leqno,twoside]{article}

\usepackage{amsmath, amsthm, amssymb,bm}

\numberwithin{equation}{section}

\newcommand{\bC}{\mathbb{C}}

\newcommand{\bZ}{\mathbb{Z}}

\newcommand{\pfa}[1]{{\textstyle \mathrm{pf}_{#1}}}
\newcommand{\dete}[1]{{\textstyle \det_{#1}}}
\newcommand{\norm}[1]{|| #1 ||}
\newcommand{\mf}[1]{\mathfrak{#1}}
\newcommand{\mr}[1]{\mathrm{#1}}
\newcommand{\mc}[1]{\mathcal{#1}}
\newcommand{\mbf}[1]{\mathbf{#1}}
\newcommand{\mbb}[1]{\mathbb{#1}}

\def\PP{\mathrm{P}}

\def\fS{\mathfrak{S}}
\def\pf{\mathrm{pf}}
\def\tr{\mathrm{tr}}

\def\fF{\mathfrak{F}}
\def\sgn{\mathrm{sgn}}

\def\per{\mathrm{per}}
\def\fo{\mathfrak{o}}
\def\fgl{\mathfrak{gl}}
\def\cB{B}
\def\bSA{{\bf SH}}

\def\cS{\mathcal{S}}

\def\({ \left( }
\def\){ \right)}
\def\[{ \left[ }
\def\]{ \right]}

\theoremstyle{plain}
\newtheorem{thm}{Theorem}[section]
\newtheorem{prop}[thm]{Proposition}
\newtheorem{lem}[thm]{Lemma}
\newtheorem{cor}[thm]{Corollary}
\theoremstyle{definition}

\newtheorem{example}{Example}[section]
\newtheorem{remark}{Remark}[section]
\theoremstyle{conjecture}
\newtheorem{conj}{Conjecture}[section]
\theoremstyle{problem}
\newtheorem{problem}[remark]{Problem}

\newcommand{\textem}[1]{{\bfseries #1}}

\title{\bfseries $\alpha$-Pfaffian, pfaffian point process \\ and shifted Schur measure}
\author{\textsc{Sho Matsumoto}}

\date{\today}

\begin{document}
\setlength{\baselineskip}{15pt}
\maketitle

\begin{abstract}
For any complex number $\alpha$ and any even-size skew-symmetric matrix $B$,
we define a generalization $\pfa{\alpha}(B)$ of the pfaffian $\pf(B)$ which we call the $\alpha$-pfaffian.
The $\alpha$-pfaffian is a pfaffian analogue of the $\alpha$-determinant studied in \cite{ST} and \cite{V}.
It gives the pfaffian at $\alpha=-1$.
We give some formulas for $\alpha$-pfaffians and study the positivity.
Further we define point processes determined by the $\alpha$-pfaffian.
Also we provide  
a linear algebraic proof of the explicit pfaffian expression obtained in \cite{Matsumoto2}
for the correlation function of the shifted Schur measure.

\par\noindent
\textem{2000 Mathematics Subject Classification} :Primary 15A15, Secondary 15A48, 05E05, 60G55
\par\noindent
\textem{Key Words} : 
$\alpha$-pfaffian, $\alpha$-determinant, pfaffian, 
pfaffian point process, Schur $Q$-function, shifted Schur measure.
\end{abstract}

%
\section{Introduction}
%

Let $\alpha$ be a complex number and $A=(a_{ij})_{1 \leq i,j \leq n}$ be a square matrix.
Then the $\alpha$-determinant of $A$ is introduced in \cite{ST} and is given by
\begin{equation} \label{alpha-determinant}
\dete{\alpha}(A)= \sum_{\sigma \in \fS_n} \alpha^{n-\nu(\sigma)} \prod_{i=1}^n a_{i \sigma(i)},
\end{equation}
where $\nu(\sigma)$ stands for the number of cycles in $\sigma \in \fS_n$ (see also \cite{V}).
Remark that in \cite{V} it is called the $\alpha$-permanent.
This gives a determinant at $\alpha=-1$ or a permanent at $\alpha=1$:
$$
\dete{-1} (A) = \sum_{\sigma \in \fS_n} \sgn(\sigma) \prod_{i=1}^n a_{i \sigma(i)}= \det (A), \qquad 
\dete{1}(A) = \sum_{\sigma \in \fS_n} \prod_{i=1}^n a_{i \sigma(i)} = \per (A).
$$

For a sequence $\bm{i}=(i_1, i_2, \dots, i_k) \in \{1,2,\dots, n\}^k$, we put
$$
A_{\bm{i}} = 
\begin{pmatrix}
a_{i_1 i_1 } & \ldots & a_{i_1 i_k} \\ \vdots & \ddots & \vdots \\ a_{i_k i_1} & \ldots & a_{i_k i_k}
\end{pmatrix}.
$$
Define a norm $\norm{A}$ of $A$ by
$$
\norm{A}= \sup_{\norm{v}=1,\ v \in \bC^n} \norm{Av},
$$
where $\norm{v}$ stands for the standard norm in $\bC^n$.
Then it is proved in \cite{V} (see also \cite{ST}) that 
if $z \in \bC$ satisfies $\norm{\alpha z A}<1$, we have
\begin{equation}\label{DetExp}
\det(I_n-\alpha z A)^{-1/\alpha} =\sum_{k=0}^\infty \frac{z^k}{k!} 
\sum_{1 \leq i_1, i_2, \dots, i_k \leq n} \dete{\alpha}(A_{i_1 i_2 \dots i_k}),
\end{equation} 
where $I_n$ is the $n$ by $n$ unit matrix.
Note that if $\alpha \in \{ -1/m ; m =1,2,\dots \}$, 
the equation \eqref{DetExp} holds without the condition for $z$
since $\det(I_n+ \frac{1}{m} z A)^m$ is a polynomial for $z$.
In particular, when $\alpha=-1$, 
the formula \eqref{DetExp} provides
the expansion of the characteristic polynomial:
\begin{equation} \label{DetExp-1}
\det(I_n+z A) = \sum_{k=0}^n z^k \sum_{1 \leq i_1 < i_2 < \dots <i_k \leq n} \det(A_{i_1 i_2 \dots i_k}).
\end{equation}

Our first aim in this paper is to obtain a pfaffian version of \eqref{DetExp}.
For that purpose,
we define a pfaffian analogue $\pfa{\alpha}(B)$ of the $\alpha$-determinant for a skew-symmetric matrix $B$
whose size is even.
This analogue $\pfa{\alpha}(B)$, which we will call the $\alpha$-pfaffian of $B$,
the $\alpha$-pfaffian gives a new expression of 
the usual pfaffian at $\alpha=-1$, i.e., $\pfa{-1}(B)=\pf(B)$ (see Section \ref{DefAlphaPf}).
Further, for a particular skew-symmetric matrix, the $\alpha$-pfaffian gives the $\alpha$-determinant.
Therefore we can regard the $\alpha$-pfaffian as a generalization of the pfaffian and the $\alpha$-determinant.

The second aim is to construct a point process determined by the $\alpha$-pfaffian.
It covers the $\alpha$-determinantal point process studied in \cite{ST} and
the usual pfaffian point process studied in \cite{BR}, \cite{R}, and \cite{S}.

The third aim is to regard the shifted Schur measure as a pfaffian point process and
to obtain a linear algebraic proof of the explicit pfaffian expression
of its correlation function 
by the same technique in \cite{BR}.
The shifted Schur measure is the measure 
which gives the probability $2^{-\ell(\lambda)} Q_{\lambda}(u) Q_{\lambda}(v)$ to
each strict partition $\lambda$,
where $\ell(\lambda)$ is the length of $\lambda$ and
$Q_{\lambda}(u)$ (resp. $Q_{\lambda}(v)$) is the Schur $Q$-function associated with $\lambda$
in the variables $u=(u_1, u_2,\dots)$ (resp. $v=(v_1, v_2, \dots)$).
Its correlation function was obtained by using vertex operators on the exterior algebra
in \cite{Matsumoto2}.
We give a new and more direct proof.

The fourth aim is to study the positivity of $\alpha$-pfaffians.
The positivity is important for existences of the point processes determined by $\alpha$-pfaffians.
In \cite{ST}, the following conjecture for $\alpha$-determinants is established.
\begin{conj}[\cite{ST}] \label{ConjST}
Let $\alpha \in [0,2] \cup \{-1/m;m=1,2,\dots\}$. Then $\dete{\alpha}(A) \geq 0$ 
whenever $A$ is a non-negative definite hermitian matrix.
\qed
\end{conj}
This is true if $\alpha \in \{2/m;m=1,2,\dots\} \cup \{0 \} \cup \{-1/m;m=1,2,\dots\}$. 
We study an $\alpha$-pfaffian version of this conjecture.

%
\section{$\alpha$-Pfaffian}
%

%
\subsection{Basic facts of the pfaffian}
%

We recall the definition of the pfaffian.
For a skew-symmetric matrix $B=(b_{kl})_{1 \leq k,l \leq 2m}$,
the pfaffian of $B$ is defined by
\begin{equation} \label{Def-pf}
\pf (B)= \sum_{\sigma \in \fF_{2m}} \sgn(\sigma) \prod_{j=1}^m b_{\sigma(2j-1) \sigma(2j)},
\end{equation}
where $\fF_{2m}$ is the subset of $\fS_{2m}$ given by
\begin{multline}
\fF_{2m}= \{\sigma=(\sigma(1),\sigma(2), \dots, \sigma(2m)) \in \fS_{2m}; \\
\sigma(2j-1)<\sigma(2j) \ (1 \leq j \leq m), \ \sigma(1)< \sigma(3)< \dots <\sigma(2m-1) \}.
\end{multline}
For convenience, we define $\pf(\emptyset)=1$ if $m=0$.
We consider a matrix of the form
$$
B=\begin{pmatrix} \cB(1,1) & \ldots & \cB(1,m) \\
\vdots & \ddots & \vdots \\ \cB(m,1) & \ldots & \cB(m,m)
\end{pmatrix},
$$
where each $\cB(r,s)$ is a $2$ by $2$ block given by
$$
\cB(r,s) = \begin{pmatrix} \cB_{00}(r,s) & \cB_{01}(r,s) \\ \cB_{10} (r,s) & \cB_{11}(r,s) \end{pmatrix}
= \begin{pmatrix} b_{2r-1,2s-1} & b_{2r-1,2s} \\ b_{2r,2s-1} & b_{2r,2s} \end{pmatrix}
\qquad \text{for $1 \leq r,s \leq m$}.
$$
We call $\cB(r,s)$ the $(r,s)$-block of $B$.
Since $b_{kl}=-b_{lk}$ for $1 \leq k,l \leq 2m$, it holds that $\cB_{ij}(r,s) = - \cB_{ji}(s,r)$
for $i,j \in \{0,1 \}$ and $1 \leq r,s \leq m$.
For a finite sequence $S=(s_1, s_2, \dots, s_k) \in \{1, \dots, m\}^k$, we put
$$
\cB[S]=
\begin{pmatrix} \cB(s_1,s_1) & \ldots & \cB(s_1,s_k) \\
\vdots & \ddots & \vdots \\ \cB(s_k,s_1) & \ldots & \cB(s_k,s_k)
\end{pmatrix}.
$$
Then $\cB[S]$ is also skew-symmetric.

Let $J_m$ be the $2m$ by $2m$ skew-symmetric matrix whose $(r,s)$-block is given by
$$
\begin{pmatrix} 0 & \delta_{rs} \\ -\delta_{rs} & 0 \end{pmatrix}.
$$
The following proposition is a pfaffian analogue of \eqref{DetExp-1}.

\begin{prop} \label{Exp-pf}
Let $B$ be any $2m$ by $2m$ skew-symmetric matrix
and $J_m$ be above.
For any $z$, we have
\begin{equation} \label{Eq-Exp-pf}
\pf(J_m+zB) = \sum_{k=0}^m z^k 
\sum_{1 \leq s_1 < s_2 < \dots < s_k \leq m } \pf(\cB[s_1, s_2, \dots, s_k]).
\end{equation}
\end{prop}

\begin{proof}
By the definition \eqref{Def-pf} of the pfaffian,
we have
\begin{align*}
\pf(J_m+zB)
=& \sum_{\sigma \in \fF_{2m}} \sgn(\sigma) 
\prod_{j=1}^m ( zb_{\sigma(2j-1) \sigma(2j)}+ (J_m)_{\sigma(2j-1) \sigma(2j)}) \\
=& \sum_{\sigma \in \fF_{2m}} \sgn(\sigma) \sum_{k=0}^m z^k 
\sum_{\begin{subarray}{c} S \subset \{1, \dots, m \} \\ \# S=k \end{subarray}}
\( \prod_{j \in S} b_{\sigma(2j-1) \sigma(2j)} \)
\( \prod_{j \in \{1, \dots ,m\}\setminus S} (J_m)_{\sigma(2j-1) \sigma(2j)} \).
\end{align*}
Since $\prod_{j \in \{1, \dots ,m\}\setminus S} (J_m)_{\sigma(2j-1) \sigma(2j)}=0$ 
unless $\sigma(2j-1)=2j-1$ and $\sigma(2j)=2j$ for all $j \in \{1, \dots, m\} \setminus S$,
we have
\begin{equation*}
\pf(J_m+zB) = \sum_{k=0}^m z^k
\sum_{1 \leq s_1 < \dots <s_k  \leq m} 
\sum_{\tau \in \fF_{2k}} \sgn(\tau) \prod_{i=1}^k b_{\tau(2s_i-1) \tau(2s_i)}
= \sum_{k=0}^m z^k 
\sum_{\begin{subarray}{c} S \subset \{1, \dots, m\}, \\ \# S=k \end{subarray}} \pf (B[S]),
\end{equation*}
where $\fF_{2k}$ acts on $\{2s_1-1, 2s_1, \dots, 2s_k-1,2s_k \}$.
\end{proof}

%
\subsection{Definition of $\alpha$-pfaffian} \label{DefAlphaPf}
%

Let $B$ be any $2m$ by $2m$ skew-symmetric matrix.
For a cycle $\tau=(k_1,k_2, \dots, k_r)$, we put
\begin{align}
& P(\cB)(\tau)= P(\cB)(k_1,k_2, \dots, k_r) \notag \\
=& \frac{1}{2} \sum_{i_1,i_2, \dots, i_r \in \bZ/2\bZ} (-1)^{i_1+i_2+\cdots+i_r} 
\cB_{i_1, i_2+1}(k_1, k_2) \cB_{i_2, i_3+1}(k_2, k_3) \cdots \cB_{i_r, i_1+1}(k_r,k_1), \label{DefinitionP}
\end{align}
where $\bZ/2\bZ=\{0,1\}$.
This is independent of the way of expressions of the cycle $\tau$:
$$
P(\cB)(k_1,k_2, \dots, k_r)=P(\cB)(k_2, \dots, k_r, k_1)= \cdots.
$$
For any complex number $\alpha$,
we define the {\it $\alpha$-pfaffian} of $B$ by 
\begin{equation} \label{Def-Pfa}
\pfa{\alpha}(B)= \sum_{\sigma \in \fS_m} \alpha^{m-\nu(\sigma)} \prod_{j=1}^{\nu(\sigma)} P(\cB)(\sigma^{(j)}),
\end{equation}
where $\sigma =\sigma^{(1)} \sigma^{(2)} \cdots \sigma^{(\nu(\sigma))}$ is the cycle decomposition of $\sigma$.

We now provide another expression of the $\alpha$-pfaffian.
Let $\tau=(k_1,k_2, \dots, k_r)$ be a cycle
and suppose that $k_1$ is the smallest among $k_1, k_2, \dots ,k_r$.
Then we put
\begin{align*}
Q(\cB)(\tau)=& Q(\cB)(k_1,k_2, \dots, k_r) \\
=&\sum_{i_2, \dots, i_r \in \bZ/2\bZ} (-1)^{i_2+\cdots+i_r} 
\cB_{0, i_2+1}(k_1, k_2) \cB_{i_2, i_3+1}(k_2, k_3) \cdots \cB_{i_r, 1}(k_r,k_1)
\end{align*}
if $r \geq 2$, or 
$Q(B)(\tau)=\cB_{01}(k_1,k_1)$ if $r=1$. 
For $\tau^{-1} = (k_1, k_r, k_{r-1}, \dots, k_2)$,
we have
$$
Q(B)(\tau^{-1}) = \sum_{i_2, \dots, i_r \in \bZ / 2\bZ} (-1)^{i_2 +\cdots+i_r}
B_{0, i_2+1}(k_1, k_r) B_{i_2,i_3+1}(k_r,k_{r-1}) \cdots B_{i_r, 1} (k_2, k_1).
$$
Since $\cB_{ij}(r,s)=-\cB_{ji}(s,r)$, changing the indices of $i_2,\dots, i_r$,
we see that
$$
= \sum_{i_2, \dots, i_r \in \bZ / 2\bZ} (-1)^{i_2 +\cdots+i_r+r}
B_{1, i_2}(k_1, k_2) B_{i_2+1,i_3}(k_2,k_3) \cdots B_{i_r+1, 0} (k_r, k_1) .
$$
Further, if we change indices as $j_p=i_p+1$, we obtain
$$
= \sum_{j_2, \dots, j_r \in \bZ / 2\bZ} (-1)^{j_2 +\cdots+j_r}
B_{1, j_2+1}(k_1, k_2) B_{j_2,j_3+1}(k_2,k_3) \cdots B_{j_r, 0} (k_r, k_1) .
$$
Therefore we get the equality
$2 P(B)(\tau) = Q(B)(\tau)+Q(B)(\tau^{-1})$.
Hence we have the

\begin{prop}
We have
\begin{equation} \label{pfa-Q}
\pfa{\alpha}(B)= \sum_{\sigma \in \fS_m} \alpha^{m-\nu(\sigma)} \prod_{j=1}^{\nu(\sigma)} Q(\cB)(\sigma^{(j)}).\qed
\end{equation}
\end{prop}

\begin{example}
We abbreviate $Q(\cB)(\tau)$ to $Q(\tau)$.

\noindent 
$\bullet$ The case where $m=1$.
We have $\pfa{\alpha}(B)= Q(1)= \cB_{01}(1,1)=b_{12}$ for any $\alpha$.

\noindent 
$\bullet$ The case where $m=2$.
Since $Q(1,2)= \sum_{i \in \bZ/2\bZ} (-1)^i \cB_{0,i+1}(1,2) \cB_{i,1}(2,1)$,
we have
\begin{align*}
\pfa{\alpha}(B)
=& Q(1)Q(2)+\alpha Q(1,2) \\
=& \cB_{01}(1,1) \cB_{01}(2,2)+ \alpha \(\cB_{01}(1,2)\cB_{01}(2,1)- \cB_{00}(1,2) \cB_{11}(2,1)\) \\
=& b_{12}b_{34} + \alpha( b_{14} b_{32}- b_{13}b_{42})
= b_{12}b_{34} + \alpha( b_{13}b_{24}-b_{14}b_{23}).\qed
\end{align*}
\end{example}

\medskip

We compare the number of terms in the expression \eqref{pfa-Q} 
with that in the definition \eqref{Def-pf} of the pfaffian.
The number of terms in \eqref{pfa-Q} as a polynomial of $\alpha$
is at most $\sum_{\sigma \in \fS_m} 2^{m-\nu(\sigma)}$.
On the other hand, the number of terms in \eqref{Def-pf} equals $(2m-1)!!$.
We can show that
\begin{equation} \label{Terms}
\sum_{\sigma \in \fS_m} 2^{m-\nu(\sigma)} = (2m-1)!!
\end{equation}
in the following Proposition \ref{PropTerms}.
Therefore the number of terms in the $\alpha$-pfaffian is equal to the number of them in the pfaffian (when $\alpha=-1$).

\begin{prop} \label{PropTerms}
For any complex number $x$, we have 
\begin{equation} \label{power}
\sum_{\sigma \in \fS_m} x^{m-\nu(\sigma)} = \prod_{k=0}^{m-1} (k x+1).
\end{equation}
In particular, $\sum_{\sigma \in \fS_m} x^{m-\nu(\sigma)} =0$ if $x = -\frac{1}{k} \ (k=1,2, \dots, m-1)$.
\end{prop}

\begin{proof}
We prove \eqref{power} by induction on $m$.
When $m=1$, both sides in \eqref{power} equal $1$.
Assume that the equality holds for $m-1$.
Put $\fS_m(m,k)= \{\sigma \in \fS_m; \sigma(m)=k \}$ for $1 \leq k \leq m$.
We define a bijective map $s$ from $\fS_m(m,k)$ to $\fS_{m-1}$ by
removing $m$ from the cycle decomposition of $\sigma \in \fS_m(m,k)$.
For example, we see that $s( (1 \ 2 \ 4)(3 \ 6 \ 5)) = (1 \ 2 \ 4)(3 \ 5)$ for $(1 \ 2 \ 4)(3 \ 6 \ 5) \in \fS_6(6,5)$.
Since $\nu(s(\sigma))$ is equal to $\nu(\sigma)$ if $\sigma \in \fS_m(m,k) \ (1 \leq k \leq m-1)$, or to
$\nu(\sigma)-1$ if $\sigma \in \fS_m(m,m)$,
we have
\begin{align*}
& \sum_{\sigma \in \fS_m} x^{m-\nu(\sigma)}= \sum_{k=1}^m \sum_{\sigma \in \fS_m(m,k)} x^{m-\nu(\sigma)}
= \sum_{k=1}^{m-1} \sum_{\pi \in \fS_{m-1}} x^{m-\nu(\pi)} + \sum_{\pi \in \fS_{m-1}}x^{m-1-\nu(\pi)} \\
=& ((m-1)x+1) \sum_{\pi \in \fS_{m-1}} x^{m-1-\nu(\pi)}.
\end{align*}
By the assumption of the induction, \eqref{power} follows immediately.
\end{proof}

The $\alpha$-pfaffian has the following invariance.

\begin{prop} \label{JBJ=B}
Let $B$ be any $2m$ by $2m$ skew-symmetric matrix.
Then we have
$$
\pfa{\alpha} ( {^t J_m B J_m}) = \pfa{\alpha}(B).
$$
\end{prop}
 
\begin{proof} 
Put $\widetilde{B}= {^tJ_m B J_m}$.
Then
it is easy to see that $\widetilde{B}_{ij} (r,s) = (-1)^{i+j} B_{i+1, j+1} (r,s)$
for $i,j=0,1$ and $r,s= 1,2, \dots, m$.
Therefore we see that 
$P(\widetilde{B}) (\tau)= P(B)(\tau)$ for any cycle $\tau$ so that
$\pfa{\alpha}(\widetilde{B}) = \pfa{\alpha}(B)$.
\end{proof}

\begin{remark}
For the usual pfaffian, the equality $\pf (^t T BT) = \pf(B) \det(T)$ holds
for any $2m$ by $2m$ skew-symmetric matrix $B$ and for any $2m$ by $2m$ matrix $T$.
But we cannot expect a similar formula for general $\alpha$-pfaffians. 
\qed
\end{remark}

%
\subsection{Expansions of the power of a pfaffian}
%

Now we give a pfaffian version of \eqref{DetExp}.

\begin{thm} \label{THM-PfExp}
Let $\alpha$ be a non-zero complex number and
$B$ be any $2m$ by $2m$ skew-symmetric matrix.
If $z \in \bC$ satisfies $\norm{\alpha z B}<1$, then we have
\begin{equation}\label{PfExp}
\pf(J_m-\alpha z B)^{-1/\alpha} =\sum_{k=0}^\infty \frac{z^k}{k!} \sum_{1 \leq s_1, s_2, \dots, s_k \leq m}
 \pfa{\alpha}(B[s_1, s_2, \dots, s_k]).
\end{equation} 
Here if $\alpha \in \{ -1/p  ; p =1,2,\dots \}$, 
the equation \eqref{PfExp} holds without the condition for $z$.
\end{thm}

To prove the theorem,
we need a lemma.

\begin{lem} \label{Lem-cycle}
For a cycle $\tau=(k_1,k_2, \dots,k_r)$, we have
$$
\tr \( (-J_1) \cB(k_1,k_2) (-J_1) \cB(k_2,k_3)\cdots (-J_1) \cB(k_r,k_1) \) = 2 P(\cB)(\tau).
$$
\end{lem}
\begin{proof}
Since
$$
-J_1 \cB(r,s)= 
\begin{pmatrix} -\cB_{10}(r,s) & -\cB_{11}(r,s) \\ \cB_{00}(r,s) & \cB_{01}(r,s) \end{pmatrix},
$$
the claim follows immediately from the definition \eqref{DefinitionP} of $P(B)(\tau)$.
\end{proof}

\begin{proof}[Proof of Theorem \ref{THM-PfExp}]
Since $\pf(B)^2 = \det(B)$, $\det(\pm J_m)=1$ and $J_m^2=-I_{2m}$,
we have
\begin{align*}
\pf (J_m -\alpha z B)^{-1/\alpha} =& \det(J_m-\alpha zB)^{-1/(2\alpha)} 
= \(\det(-J_m) \det(J_m-\alpha zB) \)^{-1/(2\alpha)} \\
=& \det(I_{2m}-\alpha z(-J_mB))^{-1/(2\alpha)}.
\end{align*}
By the assumption on $z$, we have $\norm{\alpha z J_mB}  \leq \norm{J_m} \cdot \norm{\alpha z B} <1$.
Thus we obtain
$$
\det(I_{2m}-\alpha z(-J_mB))^{-1/(2\alpha)}
= \exp \( \frac1{2\alpha} \sum_{p=1}^\infty \frac{(\alpha z)^p}{p} \tr (T^p) \), 
$$
where we put $T=-J_m B$.
If we expand the right hand side, then 
\begin{align}
\pf (J_m -\alpha z B)^{-1/\alpha} 
=& 1+ \sum_{l=1}^\infty \frac{(2\alpha)^{-l}}{l!} \sum_{p_1,p_2, \dots, p_l \geq 1} 
(\alpha z)^{p_1+ p_2+\cdots+ p_l} \frac{\tr (T^{p_1}) \tr (T^{p_2}) \cdots \tr (T^{p_l})}{p_1 p_2 \cdots p_l} \notag \\
=& 1+\sum_{k=1}^\infty (\alpha z)^k 
\sum_{l=1}^k \frac{(2\alpha)^{-l}}{l!} 
\sum_{\begin{subarray}{c} p_1,p_2, \dots, p_l \geq 1, \\ p_1+p_2+ \cdots + p_l=k \end{subarray}} 
 \frac{\tr (T^{p_1}) \tr (T^{p_2}) \cdots \tr (T^{p_l})}{p_1 p_2 \cdots p_l}. \label{eq1}
\end{align}

For a partition $\rho=(\rho_1, \rho_2, \dots, \rho_l)$ of weight $k$ and of length $\ell(\rho)=l$,
we denote by $C_\rho$ the conjugate class of cycle type $\rho$ in $\mf{S}_k$. 
It is easy to see that
$$
\frac{1}{l!}
\sum_{\begin{subarray}{c} p_1,p_2, \dots, p_l \geq 1 \\ p_1+p_2+\cdots+p_l=k  \\ 
\{ p_1, p_2, \dots, p_l \} = \{ \rho_1, \rho_2, \dots, \rho_l \} \end{subarray}} 
\frac{k!}{p_1 p_2 \cdots p_l} 
=\# C_{\rho} 
= \sum_{ \sigma \in C_{\rho} } 1,
$$
where the sum in the left hand side is all over the sequence $(p_1,p_2,\dots, p_l)$ of $l$ positive integers 
satisfying the conditions $p_1+\cdots+p_l=k$ and $\{p_1, \dots, p_l \}=\{\rho_1,\dots, \rho_l\}$. 
Hence we have
\begin{equation} \label{eq2}
\frac{1}{l!} \sum_{\begin{subarray}{c} p_1,p_2, \dots, p_l \geq 1 \\ p_1+p_2+\cdots+p_l=k \end{subarray}}
\frac{\tr(T^{p_1}) \tr(T^{p_2}) \cdots \tr(T^{p_l})}{p_1 p_2 \cdots p_l}
= \frac{1}{k!} \sum_{\begin{subarray}{c} \rho \vdash k \\ \ell(\rho)=l \end{subarray}}
\sum_{\sigma \in C_{\rho}} \tr(T^{\rho_1}) \tr(T^{\rho_2}) \cdots \tr(T^{\rho_l}).
\end{equation}
By \eqref{eq1} and \eqref{eq2},
we have
$$
\pf(J_m-\alpha z B)^{-1/\alpha} = 1+\sum_{k=1}^\infty \frac{(\alpha z)^k}{k!} 
\sum_{l=1}^k (2\alpha)^{-l} 
\sum_{\begin{subarray}{c} \rho \vdash k \\ \ell(\rho)=l \end{subarray}}
\sum_{\sigma \in C_{\rho}} \tr(T^{\rho_1}) \tr(T^{\rho_2}) \cdots \tr(T^{\rho_l}).
$$
Since 
$$
\tr(T^p)= \tr((-J_m B)^{p})= \sum_{1 \leq s_1, \dots, s_p \leq m} \tr \( (-J_1 B(s_1,s_2)) \cdots
(-J_1 B(s_p,s_1)) \)
$$
it follows from Lemma \ref{Lem-cycle} that 
\begin{align*}
\pf(J_m-\alpha z B)^{-1/\alpha} 
=& 1+\sum_{k=1}^\infty \frac{(\alpha z)^k}{k!} \sum_{1 \leq s_1, \dots, s_k \leq m}
\sum_{l=1}^k (2\alpha)^{-l} 
\sum_{\begin{subarray}{c} \rho \vdash k \\ \ell(\rho)=l \end{subarray}}
\sum_{\sigma \in C_{\rho}} 2^l \prod_{j=1}^l P(B[s_1, \dots, s_k])(\sigma^{(j)}) \\
=& 1+\sum_{k=1}^\infty \frac{z^k}{k!} \sum_{1 \leq s_1, \dots, s_k \leq m}
\sum_{l=1}^k \alpha^{k-l} 
\sum_{\begin{subarray}{c} \rho \vdash k \\ \ell(\rho)=l \end{subarray}}
\sum_{\sigma \in C_{\rho}} \prod_{j=1}^l P(\cB[s_1, \dots, s_k])(\sigma^{(j)}).
\end{align*}
If $\sigma \in C_{\rho}$, then $\nu(\sigma)=\ell(\rho)$ so that
we obtain
$$
\pf(J_m-\alpha z B)^{-1/\alpha} = 1+\sum_{k=1}^\infty \frac{z^k}{k!} \sum_{1 \leq s_1, \dots, s_k \leq m}
\pfa{\alpha} (\cB[s_1, \dots, s_k]).
$$
This completes the proof of the theorem.
\end{proof}

\begin{remark}
It is clear that 
$$
\exp \( z \sum_{s=1}^m B_{01}(s,s) \) 
= \sum_{k=0}^\infty \frac{z^k}{k!} \sum_{1 \leq s_1, \dots, s_k \leq m} \prod_{j=1}^k B_{01}(s_j, s_j)
$$
for any $z$.
We may regard this formula as a limit of \eqref{PfExp} as $\alpha \to 0$.
\qed
\end{remark}

From Proposition \ref{Exp-pf} and Theorem \ref{THM-PfExp},
we obtain a new expression of the pfaffian.

\begin{cor}
For any $2m$ by $2m$ skew-symmetric matrix $B$ we have
\begin{equation} \label{alpha-pf-pf}
\pf(B)=
\pfa{-1}(B)= \sum_{\sigma \in \fS_m} \sgn(\sigma) \prod_{j=1}^{\nu(\sigma)} Q(\cB)(\sigma^{(j)}). \qed
\end{equation}
\end{cor}

%
\subsection{Connections to $\alpha$-determinants and some formulas} \label{Connection}
%

Let $\fgl_{m}=\fgl_{m}(\bC)$ be the set of all $m$ by $m$ $\bC$-matrices 
and $\fo_{2m}=\fo_{2m}(\bC)$ be the set of all $2m$ by $2m$ skew-symmetric $\bC$-matrices.
For $A=(a_{ij})_{1 \leq i,j \leq m} \in \fgl_m$, we denote by $\omega(A) \in \fo_{2m}$ the skew-symmetric matrix 
whose $(r,s)$-block is
given by 
$$
\begin{pmatrix} 0 & a_{rs} \\ -a_{sr} & 0 \end{pmatrix}.
$$
In other words, 
\begin{equation} \label{omega}
\omega(A)=\begin{pmatrix} 0 & 1 \\ 0 & 0 \end{pmatrix} \otimes A+  \begin{pmatrix} 0 & 0 \\ -1 & 0 \end{pmatrix} 
\otimes {^tA},
\end{equation}
where $\otimes$ is the Kronecker product.
It is clear that the map $\omega:\fgl_m \to \fo_{2m}$ is linear and injective.
The following proposition assures that the $\alpha$-pfaffian is a generalization of $\alpha$-determinant.

\begin{prop} \label{pfa-dete}
For $A \in \fgl_m$,
we have $\pfa{\alpha}(\omega(A))=\dete{\alpha}(A)$.
\end{prop}

\begin{proof}
It is easy to see that
$Q(\omega(A))(\tau)=a_{k_1 k_2}a_{k_2 k_3} \cdots a_{k_r k_1}$ for a cycle $\tau=(k_1, k_2, \dots, k_r)$.
Therefore by \eqref{pfa-Q} and \eqref{alpha-determinant}
we have
$\pfa{\alpha}(\omega(A))=\dete{\alpha}(A)$.
\end{proof}

\begin{remark}
It is well known that $\pf(B)^2=\det(B)$ in the case of $\alpha=-1$.
But we cannot expect a similar formula for general $\alpha$. \qed
\end{remark}

\medskip

Theorem \ref{THM-PfExp} reduces \eqref{DetExp} by Proposition \ref{pfa-dete}.
Now we describe some formulas for $\alpha$-pfaffians and $\alpha$-determinants.

\begin{prop} \label{Alpha-Exp}
Let $B$ be any $2m$ by $2m$ skew-symmetric matrix.
Then we have
$$
\pfa{\alpha}(J_m+zB)= \sum_{S \subset \{1, \dots,m\}} z^{\# S} \pfa{\alpha}(\cB[S]).
$$
In particular, for any $m$ by $m$ matrix $A$
$$
\dete{\alpha}(I_m+zA) = \sum_{S \subset \{1, \dots,m\}} z^{\# S} \dete{\alpha}(A_S).
$$
\end{prop}

\begin{proof}
The first formula is obtained in a similar way to Proposition \ref{Exp-pf}.
The second formula is clear from the first formula and Proposition \ref{pfa-dete}.
\end{proof}

\begin{prop}
Let $B$ be any $2m$ by $2m$ skew-symmetric matrix.
Then we have
$$
\pfa{\alpha}(B)= \frac1{2^m} \sum_{i_1,\dots, i_m \in \bZ/2\bZ} 
(-1)^{i_1+\cdots+i_m} \dete{2\alpha} (\cB_{i_r, i_s+1}(r,s))_{1 \leq r,s \leq m}.
$$
\end{prop}

\begin{proof}
By the definition \eqref{Def-Pfa} of the $\alpha$-pfaffian, we have
\begin{align*}
\pfa{\alpha} (B)
=& \sum_{\sigma \in \fS_{m}} \alpha^{m-\nu(\sigma)} 2^{-\nu(\sigma)} \sum_{i_1, i_2, \dots, i_m \in \bZ/2\bZ}
(-1)^{i_1+i_2+ \cdots +i_m} \prod_{k=1}^m \cB_{i_k, i_{\sigma(k)}+1}(k, \sigma(k)) \\
=& \frac{1}{2^m} \sum_{i_1, i_2, \dots, i_m \in \bZ/2\bZ}
(-1)^{i_1+ \cdots +i_m} \sum_{\sigma \in \fS_{m}} (2\alpha)^{m-\nu(\sigma)} 
\prod_{k=1}^m \cB_{i_k, i_{\sigma(k)}+1}(k, \sigma(k)) \\
=& \frac{1}{2^m} \sum_{i_1, i_2, \dots, i_m \in \bZ/2\bZ}
(-1)^{i_1+ \cdots +i_m} \dete{2\alpha} ( \cB_{i_r,i_s+1}(r,s))_{1 \leq r,s \leq m}. 
\end{align*} 
\end{proof}

For an arbitrary matrix $T \in \mf{gl}_n$ and a non-negative integer $k$, we define $\mbb{T}^{(\alpha)}_k(T)$ by
$\mbb{T}^{(\alpha)}_k(T)= \frac{1}{2^k k!} \dete{\alpha}(Z_k(T))$, 
where 
$$
Z_k(T) = 
\begin{pmatrix}
\tr(T) & 2 & 0 & \ldots & 0 & 0 \\
\tr(T^2) & \tr(T) & 4 & \ldots & 0 & 0 \\
\vdots & \vdots & \vdots & \ddots & \vdots & \vdots \\
\tr(T^{k-2}) & \tr(T^{k-3}) & \tr(T^{k-4}) & \ldots & 2(k-2) & 0 \\
\tr(T^{k-1}) & \tr(T^{k-2}) & \tr(T^{k-3}) & \ldots & \tr(T) & 2(k-1) \\
\tr(T^k) & \tr(T^{k-1}) & \tr(T^{k-2}) & \ldots & \tr(T^2) & \tr(T) 
\end{pmatrix}
$$ 
for $k \geq 1$ and $\mbb{T}_0(T)=1$.

Then we have 

\begin{prop} \label{KW}
Let $B$ be any $2m$ by $2m$ skew-symmetric matrix.
Then we have
$$
\mbb{T}^{(\alpha)}_m(-J_m B) = \frac1{m!} \sum_{1 \leq s_1, \dots, s_m \leq m} \pfa{\alpha}(B[s_1, \dots, s_m]).
$$
In particular, $\mbb{T}^{(-1)}_m(-J_m B) = \pf (B)$.
\end{prop}

The case of $\alpha=-1$ is obtained in Lemma 2.5 in \cite{KinoshitaWakayama}.
We need a lemma to prove the proposition. 
Denote by $z_m(T)_{i j}$ the $(i,j)$-entry of the matrix $Z_m(T)$.

\begin{lem} \label{KWlemma}
Let $\rho=(\rho_1, \dots, \rho_l)$ be a partition of $m$ whose length is $l$ 
and $C_\rho$ be a conjugate class of $\mf{S}_m$ indexed by $\rho$.
For the matrix $Z_m(T)$, we have
\begin{equation} \label{KWlemmaEq}
\sum_{\sigma \in C_{\rho}} \prod_{j=1}^m z_m(T)_{j \sigma(j)} 
= 2^{m-l} \tr (T^{\rho_1}) \cdots \tr(T^{\rho_l}) \sum_{\sigma \in C_{\rho}} 1.
\end{equation}
\end{lem}

\begin{proof}
For a cycle $\tau=(1 \ 2 \ \dots \ \rho_1) (\rho_1+1 \ \dots \ \rho_1+\rho_2) \cdots$ in $C_{\rho}$,
we see that 
\begin{align*}
\prod_{j=1}^m z_m(T)_{j \tau(j)} 
=& 2 \cdot 4 \cdots 2(\rho_1-1) \cdot \tr(T^{\rho_1}) \cdot 2(\rho_1+1) \cdots
2(\rho_1+\rho_2-1) \tr(T^{\rho_2}) \cdots \\
=& 2^{m-l} \frac{m!}{\rho_1 (\rho_1+\rho_2) \cdots (\rho_1+ \rho_2+ \cdots +\rho_l)} 
\tr (T^{\rho_1}) \cdots \tr(T^{\rho_l}).
\end{align*}
Since $z_m(T)_{ij}=0$ if $j>i+1$, 
only the terms for the permutations whose forms are $\sigma= (1 \ 2 \ \dots \ a_1) (a_1+1 \ \dots \ a_1+a_2) \cdots$
remain in the left hand side of \eqref{KWlemmaEq},
where $a_1, a_2, \dots$ are permutations for $\rho_1, \rho_2, \dots$.
Therefore
\begin{align*}
& \sum_{\sigma \in C_{\rho}} \prod_{j=1}^m z_m(T)_{j \sigma(j)} \\
=& 2^{m-l} m! \tr(T^{\rho_1}) \cdots \tr(T^{\rho_l}) 
\( \prod_{r \geq 1} \frac{1}{m_r(\rho)!} \)
\sum_{\pi \in \mf{S}_l} 
\frac{1}{\rho_{\pi(1)}( \rho_{\pi(1)}+\rho_{\pi(2)}) \cdots (\rho_{\pi(1)} +\rho_{\pi(2)} + \cdots +\rho_{\pi(l)})},
\end{align*}
where $m_r(\rho)$ is the multiplicity of $r$ in $\rho=(\rho_1, \dots, \rho_l)$.
It is easy to see that this equals
$$
= 2^{m-l} m! \tr(T^{\rho_1}) \cdots \tr(T^{\rho_l}) 
\( \prod_{r \geq 1} \frac{1}{m_r(\rho)!}\) \frac{1}{\rho_1 \cdots \rho_l} 
= 2^{m-l} \tr(T^{\rho_1}) \cdots \tr(T^{\rho_l}) \# C_{\rho}.
$$
We have proved the lemma.
\end{proof}

\begin{proof}[Proof of Proposition \ref{KW}]
We have
$$
\mbb{T}^{(\alpha)}_m(-J_m B)
= \frac{1}{2^m m!} \sum_{\rho \vdash m} \alpha^{m-\ell(\rho)}
\sum_{\sigma \in C_{\rho}} \prod_{j=1}^m z_m(-J_m B)_{j \sigma(j)}, 
$$
where $\rho \vdash m$ means that $\rho$ is a partition of $m$.
By Lemma \ref{KWlemma}, 
$$
= \frac{1}{m!} \sum_{\rho \vdash m} \alpha^{m-\ell(\rho)} 2^{-\ell(\rho)}
\tr((-J_m B)^{\rho_1}) \cdots \tr((-J_m B)^{\rho_l}) \sum_{\sigma \in C_\rho}1.
$$
In a similar way to the proof of Theorem \ref{THM-PfExp},
we can obtain the first claim in the proposition.

In the case of $\alpha=-1$, 
we see that
$$
\mbb{T}^{(-1)}_m(-J_m B)= \frac{1}{m!} 
\sum_{\begin{subarray}{c} 1 \leq s_1, \dots, s_m \leq m, \\ \text{$s_1, \dots, s_m$ are distinct} \end{subarray}}
 \pf(B[s_1, \dots, s_m])
= \pf(B).
$$
\end{proof}

\medskip

\begin{remark}
Let $\lambda$ be a partition of $n$ and $A=(a_{ij})$ be a matrix in $\mf{gl}_n$.
Then the immanant of $A$ is defined by 
$ \mr{Imm}_\lambda (A)=\sum_{\sigma \in \mf{S}_n} \chi^{\lambda}(\sigma) \prod_{i=1}^n a_{i \sigma(i)}$,
where $\chi^{\lambda}$ is the irreducible character of $\mf{S}_n$ associated with $\lambda$,
see \cite{J}.
It is easy to see that $\mr{Imm}_{(n)} (A) = \per (A)$ and $\mr{Imm}_{(1^n)} (A)= \det (A)$.
We can define a pfaffian analogue of the immanant by
$$
\pf^{\lambda} (B) = \sum_{\sigma \in \mf{S}_m} \chi^{\lambda} (\sigma) \prod_{j=1}^{\nu(\sigma)} Q(B)(\sigma^{(j)}),
$$
where $B \in \mf{o}_{2m}$ and $\lambda$ is a partition of $m$.
We immediately see $\pf^{(m)}(B)= \pfa{1} (B)$ and $\pf^{(1^m)}(B) = \pfa{-1}(B) = \pf (B)$.
As like Proposition \ref{pfa-dete}, we obtain 
$\pf^{\lambda}(\omega(A)) = \mr{Imm}_\lambda (A)$ for $A \in \mf{gl}_m$ and a partition $\lambda$ of $m$,
that is, $\pf^{\lambda}$ is an extension of the immanant.
\qed
\end{remark}

%
\section{$\alpha$-Pfaffian point processes} \label{PPP}
%

Let $\mf{X}$ be a countable set.
Let $L$ be a map 
$$
L: \mf{X} \times \mf{X} \rightarrow \fgl_2 (\bC); (x,y) 
\mapsto L(x,y)= \begin{pmatrix} L_{00}(x,y) & L_{01}(x,y) \\ L_{10}(x,y) & L_{11}(x,y) \end{pmatrix}
$$
such that $L_{ij} (x,y) = -L_{ji}(y,x)$ for any $i,j \in \{0,1\}$ and $x,y \in \mf{X}$. 
Such $L$ is called a skew-symmetric matrix kernel on $\mf{X}$, see \cite{R,S}.
We regard the map $L$ as an operator on the Hilbert space $\ell^2(\mf{X}) \oplus \ell^2(\mf{X})$.
Then $L$ is a matrix whose blocks are indexed by elements in $\mf{X} \times \mf{X}$.
For $n$-point $x_1,x_2,\dots, x_n \in \mf{X}$, we denote by $L[x_1,x_2, \dots, x_n]$ the $2n$ by $2n$ skew-symmetric
matrix $(L(x_i,x_j))_{1 \leq i,j \leq n}$. 

Let $J$ be the skew-symmetric matrix kernel determined by $J(x,y)= J_1 \delta_{x,y}$ and $\alpha$ be a real number.
Suppose that $\pfa{\alpha}(L[x_1, \dots, x_n])$
are non-negative for all $(x_1,x_2, \dots, x_n) \in \mf{X}^n \ ( n \geq 1)$ and that $\norm{\alpha L} <1$,
where $\norm{\cdot}$ is the operator norm with respect to the Hilbert space $\ell^2(\mf{X}) \oplus \ell^2(\mf{X})$.
By Theorem \ref{THM-PfExp}, we have
\begin{equation} \label{PfExpInfinite}
\pf(J-\alpha L)^{-1/\alpha} = \sum_{n=0}^\infty \frac{1}{n!} 
\sum_{x_1, \dots, x_n \in \mf{X}} \pfa{\alpha}(L[x_1, \dots, x_n]).
\end{equation}
Here if $\# \mf{X} =\infty$ then we put 
$\pf(J+L) := \sum_{\begin{subarray}{c} X \subset \mf{X}, \\ \# X < \infty \end{subarray}} \pf (L[X])$.
We consider a probability density on $\bigcup_{n=0}^\infty \mf{X}^n$ defined by
\begin{equation} 
p^{(\alpha)}_{L} (x_1, \dots, x_n)
 = \frac{1}{n!} \pf(J-\alpha L)^{1/\alpha} \pfa{\alpha}(L[x_1, \dots, x_n]) \qquad \text{on $\mf{X}^n$}.
\end{equation}
By \eqref{PfExpInfinite}, we have 
$\sum_{n=0}^\infty \sum_{x_1, \dots ,x_n \in \mf{X}} p^{(\alpha)}_L(x_1, \dots, x_n)=1$.

Put $\mf{Q}(\mf{X})=\bigcup_{n=0}^\infty \mf{X}^n /\sim$,
where we write $(x_1, \dots, x_n) \sim (y_1, \dots, y_n)$ for two elements in $\mf{X}^n$
if and only if there exists a permutation $\sigma \in \mf{S}_n$ such that $y_j=x_{\sigma(j)} \ (1 \leq j \leq n)$.
We define a point process $\pi^{(\alpha)}_L$ on $\mf{X}$ by the equality
\begin{equation}
\sum_{X \in \mf{Q}(\mf{X})} \pi^{(\alpha)}_L(X) f(X)
= \sum_{n=0}^\infty \sum_{x_1, \dots, x_n \in \mf{X}} p^{(\alpha)}_L(x_1, \dots, x_n) f(x_1) \cdots  f(x_n)
\end{equation}
for any test function $f$ on $\mf{X}$.
Here we write $f(X)=f(x_1) \cdots f(x_n)$ for any sequence $(x_1, \dots, x_n)$ in a class $X$.
More explicitly,
$\pi^{(\alpha)}_L$ is a probability measure on $\mf{Q}(\mf{X})$ defined by
\begin{equation} \label{PointProcess}
\pi^{(\alpha)}_L(X) = \frac1{\prod_{x \in \mf{X}} m_x(X)!} \pf(J-\alpha L)^{1/\alpha} \pfa{\alpha}(L[X]) \qquad
\text{for $X \in \mf{Q}(\mf{X})$},
\end{equation}
where $m_x(X)$ is the multiplicity of $x$ in $X$.
We call this measure $\pi^{(\alpha)}_L$ the {\it $\alpha$-pfaffian point process} on $\mf{X}$
determined by a skew-symmetric matrix kernel $L$.

The $n$-point correlation function is defined by
\begin{equation}
\rho^{(\alpha)}_L(x_1, \dots, x_n)= \sum_{k=0}^\infty \frac{(n+k)!}{k!} \sum_{x_{n+1}, \dots, x_{n+k} \in \mf{X}} 
p^{(\alpha)}_L(x_1, \dots, x_n,x_{n+1}, \dots ,x_{n+k})
\end{equation}
for $(x_1, \dots, x_n) \in \mf{X}^n$.

\begin{thm} \label{THM-Correlation}
Let $L$, $p^{(\alpha)}_L$, and $\rho^{(\alpha)}_L$ be as above.
Put $K_\alpha = -\alpha^{-1} (J + (J-\alpha L)^{-1} )$ and
$\widetilde{K}_{\alpha}=L(I+\alpha JL)^{-1}$.
Assume that $\pfa{\alpha}(K_\alpha[x_1, \dots, x_n])$ is non-negative for any  $x_1, \dots, x_n \in \mf{X}$.
Then we have $\rho^{(\alpha)}_L(x_1, \dots, x_n) = \pfa{\alpha}(K_\alpha[x_1, \dots ,x_n])
=\pfa{\alpha}(\widetilde{K}_\alpha[x_1, \dots ,x_n])$.
\end{thm}

\begin{proof}
Put $K'_\alpha = -J \widetilde{K}_\alpha J$.
Then it follows that
$(J+\alpha K'_\alpha)(J- \alpha L) = -I -\alpha JL + \alpha (-J \widetilde{K}_\alpha J) J(I+\alpha JL)
= -I$ so that
$K'_\alpha = -\alpha^{-1} ( J+ (J-\alpha L)^{-1} ) =K_\alpha$.
Therefore by Proposition \ref{JBJ=B},
it is enough to prove 
$\rho^{(\alpha)}_L(x_1, \dots, x_n) = \pfa{\alpha}(\widetilde{K}_\alpha[x_1, \dots ,x_n])$.

We abbreviate $p^{(\alpha)}_L$ to $p$ simply.
Let $f(x)$ be a function on $\mf{X}$ such that $f_0(x) =f(x)-1$ is finitely supported.
We have
\begin{align*}
 & \sum_{m=0}^\infty \sum_{y_1, \dots ,y_m \in \mf{X}} p(y_1, \dots, y_m) \prod_{j=1}^m f(y_j) 
=  \sum_{m=0}^\infty \sum_{y_1, \dots ,y_m \in \mf{X}} p(y_1, \dots, y_m) \prod_{j=1}^m (1+f_0(y_j)) \\
=& \sum_{m=0}^\infty \sum_{n=0}^m \begin{pmatrix} m \\ n \end{pmatrix} \sum_{x_1, \dots, x_n \in \mf{X}} 
\prod_{j=1}^n f_0(x_j) \sum_{x_{n+1}, \dots, x_m \in \mf{X}} p(x_1, \dots, x_n, x_{n+1}, \dots, x_m) \\
=& \sum_{n=0}^\infty \sum_{k=0}^\infty \frac{(n+k)!}{n! k!}  \sum_{x_1, \dots, x_n \in \mf{X}} 
\prod_{j=1}^n f_0(x_j) \sum_{x_{n+1}, \dots, x_{n+k} \in \mf{X}} p(x_1, \dots, x_n, x_{n+1}, \dots, x_{n+k}) \\
=& \sum_{n=0}^\infty \frac{1}{n!} \sum_{x_1, \dots, x_n \in \mf{X}} 
\rho^{(\alpha)}_L (x_1, \dots, x_n) \prod_{j=1}^n f_0(x_j).
\end{align*}
When we put $\sqrt{f} (x) = \sqrt{f(x)}$ and identify $\sqrt{f}$ with the matrix 
with $(x,y)$-block $\sqrt{f}(x) I_2 \delta_{x,y}$, we have $\det(\sqrt{f}[x_1, \dots, x_n])= f(x_1) \cdots f(x_n)$.
Therefore we get
\begin{align*}
& \sum_{m=0}^\infty \sum_{y_1, \dots ,y_m \in \mf{X}} p(y_1, \dots, y_m) \prod_{j=1}^m f(y_j) \\
=& 
\sum_{m=0}^\infty \frac1{m!} \sum_{y_1, \dots ,y_m \in \mf{X}} \pfa{\alpha} ((\sqrt{f}L \sqrt{f})[y_1, \dots, y_m])
\pf (J-\alpha L)^{1/\alpha}  \\
=& \pf(J-\alpha \sqrt{f} L \sqrt{f})^{-1/\alpha} \pf(J-\alpha L)^{1/\alpha}
\end{align*}
by \eqref{PfExpInfinite}.
Further we see that 
\begin{align*}
& \det(J-\alpha \sqrt{f} L \sqrt{f}) \det( J-\alpha L)^{-1}
= \det( (J-\alpha f L)(J-\alpha L)^{-1})  \\
=& \det((-J)(J-\alpha f L)(J-\alpha L)^{-1} J) = \det( (I+\alpha (1+f_0) JL) (I+\alpha J L)^{-1} ) \\
=& \det(I +f_0 J \widetilde{K}_\alpha) =\det( J-\alpha f_0 \widetilde{K}_\alpha) 
= \det(J-\alpha \sqrt{f_0} \widetilde{K}_\alpha \sqrt{f_0}).
\end{align*}
Hence, by \eqref{PfExpInfinite} again,
\begin{align*}
& \sum_{m=0}^\infty \sum_{y_1, \dots ,y_m \in \mf{X}} p(y_1, \dots, y_m) \prod_{j=1}^m f(y_j) 
=  \pf(J-\alpha \sqrt{f_0} L \sqrt{f_0})^{-1/\alpha} \pf(J-\alpha L)^{1/\alpha} \\
=& \pf(J-\alpha \sqrt{f_0} \widetilde{K}_\alpha \sqrt{f_0})^{-1/\alpha} 
= \sum_{n=0}^\infty \frac{1}{n!} \sum_{x_1, \dots, x_n \in \mf{X}} 
\pfa{\alpha} (\widetilde{K}_\alpha[x_1,\dots, x_n]) \prod_{j=1}^n f_0(x_j).
\end{align*}
Finally we have
$$
\sum_{n=0}^\infty \frac{1}{n!}\sum_{x_1, \dots, x_n \in \mf{X}}  
\rho^{(\alpha)}_L (x_1, \dots, x_n) \prod_{j=1}^n f_0(x_j)
= \sum_{n=0}^\infty \frac{1}{n!} \sum_{x_1, \dots, x_n \in \mf{X}} 
\pfa{\alpha} (\widetilde{K}_\alpha[x_1,\dots, x_n]) \prod_{j=1}^n f_0(x_j)
$$
for any finitely supported function $f_0$.
This implies the theorem.
\end{proof}

In general,
we call the matrix $K$ such that $\rho(x_1, \dots, x_n) = \pf(K[x_1, \dots, x_n])$
the correlation kernel for a correlation function $\rho$.
The correlation kernel is not uniquely determined.

The $\alpha$-determinant version of the theorem above is obtained in \cite{ST}.

\begin{example}
Let $A$ be a trace class operator on $\mf{X}$.
Define a skew-symmetric matrix kernel $L$ by 
$$
L(x,y)= \begin{pmatrix} 0 & A(x,y) \\ -A(y,x) & 0 \end{pmatrix}
$$
for $x,y \in \mf{X}$.
Then $\pfa{\alpha}(L[X]) = \dete{\alpha} (A_X)$ by Proposition \ref{pfa-dete}, where $A_X= (A(x,y))_{x,y \in X}$.
Therefore the probability measure $\pi^{(\alpha)}_L$ is nothing but
the $\alpha$-determinantal point process determined by $A$, which studied in \cite{ST}.
The correlation function is given as 
$\rho^{(\alpha)}_L(x_1, \dots, x_n) = 
\dete{\alpha}(B_{\alpha}(x_i, x_j))_{1 \leq i,j \leq n}$,
where $B_\alpha = A(I-\alpha A)^{-1} = -\alpha^{-1} (I -(I+A)^{-1})$. 
\qed
\end{example}

\begin{remark}
If 
$\pfa{\alpha}(L[x_1, \dots, x_n])$ is non-negative for $x_1, \dots, x_n \in \mf{X}$ 
and
$\pfa{\alpha} (J+L)$ is non-negative, then
we can also define the point process on $\mf{P}(\mf{X})=\{ X \subset \mf{X} ; \text{ $\#X$ is finite} \}$
from Proposition \ref{Alpha-Exp}
by
$\Pi^{(\alpha)}_L (X) =  \pfa{\alpha} (J+L)^{-1} \pfa{\alpha}(L[X])$
for $X \in \mf{P}(\mf{X})$.
\qed
\end{remark}

%
%
\section{Shifted Schur measure and pfaffian point process}
%
%

In this section,
we obtain a linear algebraic proof of the explicit pfaffian expression of 
the correlation function for the shifted Schur measure.
It is obtained in \cite{Matsumoto2} via the vertex operators on the exterior algebra.

%
\subsection{Definition of the shifted Schur measure}
%

Let $u=(u_1,u_2, \dots)$ and $v=(v_1, v_2, \dots)$ be variables.
Let $\mc{D}$ be the set of all strict partitions and $\ell(\lambda)$ 
be the length of a partition $\lambda \in \mc{D}$ (see \cite{Mac}).
Put
$$
Q_u(z)=  \prod_{j=1}^\infty \frac{1+u_j z}{1-u_j z}.
$$
The Schur $Q$-function $Q_{\lambda}(u)$ associated with $\lambda \in \mc{D}$ is defined as a coefficient of 
$z_1^{\lambda_1} \cdots z_n^{\lambda_n}$
in the formal series expansion of
$$
Q_u(z_1) \cdots Q_u(z_n) \prod_{1 \leq i< j \leq n} \frac{z_i-z_j}{z_i+z_j},
$$
where $n \geq \ell(\lambda)$ and
$\frac{z-w}{z+w} = 1+ 2 \sum_{k=1}^\infty (-1)^k z^{-k}w^k$.

The shifted Schur measure (see \cite{Matsumoto2, TW2004}) is the probability measure on $\mc{D}$ given by
$$
\PP_{\mr{SS}} (\lambda) = \frac{1}{Z_{\mr{SS}}} 2^{-\ell(\lambda)} Q_{\lambda}(u) Q_{\lambda}(v),
$$
where $Z_{\mr{SS}}$ is the constant determined by the Cauchy identity (\cite[III-8]{Mac})
$$
Z_{\mr{SS}}= \sum_{\lambda \in \mc{D}} 2^{-\ell(\lambda)} Q_{\lambda}(u) Q_{\lambda}(v) 
= \prod_{i,j=1}^\infty \frac{1+u_i v_j}{1-u_iv_j}.
$$
The correlation function for the shifted Schur measure is defined by 
$$
\rho_{\mr{SS}}(\lambda) =  \PP_{\mr{SS}} 
\( \left\{ \mu \in \mc{D} ;  \ \mu \supset \lambda \right\} \)
= \sum_{ \mu \in \mc{D},\ \mu \supset \lambda} \PP_{\mr{SS}}(\mu)
\qquad
\text{for $\lambda \in \mc{D}$},
$$
where $\mu \supset \lambda$ means
$\{ \mu_1, \dots, \mu_{\ell(\mu)} \} \supset \{\lambda_1, \dots, \lambda_{\ell(\lambda)} \}$.
The aim in this section is to obtain a new proof of the following theorem in \cite{Matsumoto2}.

\begin{thm} \label{THM-SS}
For any strict partition $\lambda=(\lambda_1, \dots, \lambda_l) \in \mc{D}$,
we have
$$
\rho_{\mr{SS}}(\lambda ) = \pf (\mc{K}[\lambda_1, \dots, \lambda_l]),
$$
where $\mc{K}$ is a skew-symmetric matrix kernel whose blocks are
$$
\mc{K}(r,s) 
= \begin{pmatrix} 
\mc{K}_{00} (r,s) &\mc{K}_{01} (r,s) \\ \mc{K}_{10} (r,s) & \mc{K}_{11} (r,s)
\end{pmatrix}
\qquad \text{for $r,s \geq 1$}.
$$
The each entry is given as follows:
$$
\mc{K}_{00}(r,s) = \frac{1}{2} [z^r w^s] \frac{Q_u(z) Q_u(w)}{Q_v(z^{-1}) Q_v(w^{-1})} \frac{z-w}{z+w},
$$
where $\frac{z-w}{z+w} = 1 + 2\sum_{k=1}^\infty (-1)^k z^{-k} w^{k}$
and $[z^r w^s]$ stands for the coefficient of $z^r w^s$,
and 
$$
\mc{K}_{01}(r,s) = -\mc{K}_{10}(s,r) =
\frac{1}{2} [z^r w^s] \frac{Q_u(z) Q_v(w)}{Q_v(z^{-1}) Q_u(w^{-1})} \frac{zw+1}{zw-1}, 
$$
where 
$\frac{zw+1}{zw-1} =
 -\(1 + 2 \sum_{k=1}^\infty  z^{k} w^{k}\)$.
Finally,
$$
\mc{K}_{11} (r,s) = \frac{1}{2} [z^r w^s] \frac{Q_v(z) Q_v(w)}{Q_u(z^{-1}) Q_u(w^{-1})} \frac{w-z}{w+z},
$$
where $\frac{w-z}{w+z} = -\( 1 + 2\sum_{k=1}^\infty (-1)^k z^{-k} w^{k} \)$.
\end{thm}

We regard the shifted Schur measure as a pfaffian point process.
For that purpose, 
we recall the pfaffian point process,
which is the case of $\alpha=-1$ in the $\alpha$-pfaffian point process defined in Section \ref{PPP}.

Let $\mf{P}(\mf{X})=\{ X \subset \mf{X} ; \text{ $\#X$ is finite} \}$.
We define the pfaffian point process $\pi_L$ on $\mf{X}$ determined 
by a skew-symmetric matrix kernel $L$ as the measure on $\mf{P}(\mf{X})$
given by
$$
\pi(X)=\pi_L(X)= \frac{\pf(L[X])}{\pf(J+L)} \qquad \text{for $X \in \mf{P}(\mf{X})$}.
$$
Then by Theorem \ref{THM-Correlation}
its correlation function is 
\begin{equation} \label{Pfaffian-correlation}
\rho(X)= \sum_{Y \in \mf{P}(\mf{X}), \ Y \supset X} \pi(Y) = \pf (K[X]),
\end{equation}
where $K=J+(J+L)^{-1}$.

More generally, let $\mf{Y}$ be a subset in $\mf{X}$ such that $\mf{Y}^c = \mf{X} \setminus \mf{Y}$ is finite.
Then we can define a point process on $\mf{Y}$ by
\begin{equation} \label{CPPP}
\pi_{L, \mf{Y}} (X) = \frac{\pf(L[X \cup \mf{Y}^c])}{\pf (J[\mf{Y}]+L)} 
\qquad \text{for $X \in \mf{P}(\mf{Y})$}.
\end{equation}
Here we identify $J[\mf{Y}]$ with the block matrix
$\begin{pmatrix} J & 0 \\ 0 & 0 \end{pmatrix}$, 
where the blocks correspond to the partition $\mf{X} = \mf{Y} \sqcup \mf{Y}^c$.
This process is called the conditional pfaffian point process determined by $L$.
The correlation function is given as follows.

\begin{prop}[\cite{BR}] \label{ConditionalPPP}
Let $\rho_{L,\mf{Y}}$ be the correlation function on $\mf{Y}$ determined by a skew-symmetric matrix kernel $L$.
Then we have
$\rho_{L,\mf{Y}} (X) = \pf (K[X])$ for $X \in \mf{P}(\mf{Y})$,
where 
$$
K=J[\mf{Y}] +(J[\mf{Y}]+L)^{-1} \Big|_{\mf{Y} \times \mf{Y}}. \qed
$$
\end{prop}

%
\subsection{Proof of Theorem \ref{THM-SS}: 1st step} 
%

Let $u=(u_1, \dots, u_n)$ and $v=(v_1, \dots, v_n)$, 
where $n$ is even.
The case that $n$ is odd is similar.
We define a bijective map $\phi$ from $\mc{D}$ to 
$\mf{P}^{\mr{even}}(\bZ_{\geq 0})= \{ X \in \mf{P}(\bZ_{\geq 0}) \ ; \ \text{$\# X$ is even} \}$ by
$$
\phi(\lambda) = 
\begin{cases} \{ \lambda_1, \dots, \lambda_{\ell(\lambda)} \}, & \text{if $\ell(\lambda)$ is even}, \\
\{ \lambda_1, \dots, \lambda_{\ell(\lambda)},0 \}, & \text{if $\ell(\lambda)$ is odd}.
\end{cases}
$$
First we prove the following proposition.

\begin{prop} \label{Prop-L}
Define a skew-symmetric matrix kernel $L$ on 
$\mf{X} = \{1,2, \dots, n \} \sqcup \mf{Y}$
by
$$
L = \begin{pmatrix} \mc{V} & \mc{W} \eta^{-\frac{1}{2}} \\ -\eta^{-\frac{1}{2}} {^t \mc{W}} & O \end{pmatrix}, 
$$
where $\mf{Y}=\bZ_{\geq 0}$, $\mc{V}= (\mc{V}(i,j))_{1 \leq i,j \leq n}$ 
and $\mc{W} =(\mc{W}(i,r))_{1 \leq i \leq n,\ r \in \bZ_{\geq 0}}$.
Their entries are given by
$$
\mc{V} (i,j) = \begin{pmatrix} - \frac{u_i-u_j}{u_i+u_j}& 0 \\ 0 & \frac{v_i-v_j}{v_i+v_j} \end{pmatrix}, \qquad
\mc{W} (i,r) = \begin{pmatrix} -u^r_i & 0 \\ 0 & v^r_i \end{pmatrix}.
$$
Further $\eta$ is the matrix whose block is given by 
$$
\eta (r,s) = \delta_{rs} \begin{pmatrix} \eta(r) & 0 \\ 0 & \eta(r) \end{pmatrix}
\qquad \text{for $r,s \in \bZ_\geq 0$},
$$
where $\eta(r)$ is equal to $1$ if $r=0$, or to $\frac{1}{2}$ if $r \geq 1$.
Then for $X \in \mf{P} ( \bZ_{\geq 0})$ we have
$$
\pf (L[X \cup \{1,2,\dots, n \}])
=2^{-\ell(\lambda)} Q_{\lambda} (u) Q_{\lambda}(v) 
\prod_{1 \leq i < j \leq n} \( \frac{u_i-u_j}{u_i+u_j} \)   \prod_{1 \leq i < j \leq n} \( \frac{v_i-v_j}{v_i+v_j} \)
$$ 
if $\# X$ is even and $X=\phi(\lambda)$,
and $\pf (L[X \cup \{1,2,\dots, n \}])=0$ if $\# X$ is odd.
\end{prop}

In order to prove the proposition,
we need the lemma (see \cite[III-8]{Mac}).

\begin{lem}[\cite{N}] \label{Nimmo}
The Schur Q-function can be written as 
$$
Q_{\lambda}(u) = \pm 2^{\ell(\lambda)} \pf  (A_u(\omega(\lambda)))\cdot
\prod_{1 \leq i < j \leq n} \( \frac{u_i-u_j}{u_i+u_j} \)^{-1},
$$
where
the sign depends on only $\lambda$
and we set 
$$
A_u(X)= \begin{pmatrix} V(u) & W(u)_X \\ - ^t W(u)_X & O_X \end{pmatrix}
$$
for $X \in \mf{P}(\bZ_{\geq 0})$.
Here 
$V(u)= \( \frac{u_i-u_j}{u_i+u_j} \)_{1 \leq i,j \leq n}$
and $W(u)_X = (u_i^{x})_{1 \leq i \leq n, x \in X}$.
\qed
\end{lem}

\begin{proof}[Proof of Proposition \ref{Prop-L}]
Put $\mc{W}[X]=(\mc{W}(i,x))_{1 \leq i \leq n, x \in X}$ for $X \in \mf{P}(\bZ_{\geq 0})$.
Since 
\begin{align*}
 & L[X \cup \{1,2, \dots, n\}]
= \begin{pmatrix} \mc{V} & \mc{W}[X] \  \eta^{-\frac{1}{2}}[X] \\ 
-\eta^{-\frac{1}{2}}[X] \ {^t \mc{W}[X]} & O_X \end{pmatrix} \\
=& \begin{pmatrix} I_{2n} & O \\ O & \eta^{-\frac{1}{2}} [X] \end{pmatrix} 
\begin{pmatrix} \mc{V} & \mc{W}[X]  \\ -{^t \mc{W}[X]} & O_X \end{pmatrix}
\begin{pmatrix} I_{2n} & O \\ O & \eta^{-\frac{1}{2}} [X] \end{pmatrix}
\end{align*}
for $X \in \mf{P}(\bZ_{\geq 0})$,
we have 
\begin{align*}
\pf (L[X \cup \{1,2, \dots, n\}]) 
=& \pf  \begin{pmatrix} \mc{V} & \mc{W}[X]  \\ -{^t \mc{W}[X]} & O_X \end{pmatrix}  
\det  \begin{pmatrix} I_{2n} & O \\ O & \eta^{-\frac{1}{2}} [X] \end{pmatrix} \\
=&  \prod_{x \in X} \eta(x)^{-1} \cdot \pf  \begin{pmatrix} \mc{V} & \mc{W}[X]  \\ -{^t \mc{W}[X]} & O_X \end{pmatrix}.
\end{align*}
Changing orders of rows and columns,
we have 
\begin{align*}
\pf  \begin{pmatrix} \mc{V} & \mc{W}[X]  \\ -{^t \mc{W}[X]} & O_X \end{pmatrix}
=& (-1)^{(n+\# X)(n+\# X-1)/2} \pf \begin{pmatrix} -A_u(X) & O \\ O & A_v(X) \end{pmatrix} \\
=& (-1)^{(n+\# X)(n+\# X-1)/2+ (n+\# X)/2} \pf (A_u(X)) \pf (A_v(X)).
\end{align*}
Here if $\# X$ is odd then $\pf (A_u(X))=0$. 
Hence we may assume $\# X$ is even.
Put $\lambda = \phi^{-1} (X)$. 
Then it follows from Lemma \ref{Nimmo} that
\begin{align*}
 &\pf  \begin{pmatrix} \mc{V} & \mc{W}[X]  \\ -{^t \mc{W}[X]} & O_X \end{pmatrix}
= \pf (A_u(X)) \pf (A_v(X)) \\
=& (2^{-\ell(\lambda)})^2  Q_{\lambda}(u) Q_{\lambda}(v) 
\prod_{1 \leq i < j \leq n} \( \frac{u_i-u_j}{u_i+u_j} \)   \prod_{1 \leq i < j \leq n} \( \frac{v_i-v_j}{v_i+v_j} \).
\end{align*}
Since $\prod_{x \in X} \eta(x) = 2^{-\ell(\lambda)}$ for $X=\phi(\lambda)$,
we get the proposition.
\end{proof}

%
\subsection{Proof of Theorem \ref{THM-SS}: 2nd step} 
%

The shifted Schur measure is nothing but
the conditional pfaffian point process determined by $L$ in Proposition \ref{Prop-L}.
Thus from Proposition \ref{ConditionalPPP}
we have to obtain an explicit expression
of $K=J[\mf{Y}] + (J[\mf{Y}] + L)^{-1} \Big|_{\mf{Y} \times \mf{Y}}$.

\begin{lem}
We have
$$
\begin{pmatrix} A & B \\ C & D \end{pmatrix}^{-1} =
\begin{pmatrix} -\mc{M}^{-1} & \mc{M}^{-1} B D^{-1} \\
D^{-1} C \mc{M}^{-1} & D^{-1} - D^{-1} C \mc{M}^{-1} B D^{-1} \end{pmatrix},
$$
where $\mc{M}= B D^{-1} C -A$,
if $D$ and $\mc{M}$ are invertible.
\end{lem}

\begin{proof}
This is straightforward, see e.g. \cite{BR}.
\end{proof}

By this lemma, the kernel
$K=J[\mf{Y}] + (J[\mf{Y}] + L)^{-1} \Big|_{\mf{Y} \times \mf{Y}}$
is equal to
$K = J[\mf{Y}] \eta^{-\frac{1}{2}} {^t \mc{W}} \mc{M}^{-1} \mc{W} \eta^{-\frac{1}{2}} J[\mf{Y}]$,
where $\mc{M}= \mc{W} \eta^{-\frac{1}{2}} J[\mf{Y}] \eta^{-\frac{1}{2}} {^t \mc{W}} - \mc{V}$.
We may replace $K$ by 
$$
-\eta^{-\frac{1}{2}} {^t \mc{W}} \mc{M}^{-1} \mc{W} \eta^{-\frac{1}{2}}
$$
from Proposition \ref{JBJ=B}.

\begin{prop} \label{Prop-MI}
Write the skew-symmetric matrix kernel $\mc{M}^{-1}$ on $\{1, 2, \dots, n \}$
as
$$
\mc{M}^{-1} (k,l) = \begin{pmatrix} \mc{M}^{-1}_{00}(k,l) & \mc{M}^{-1}_{01}(k,l) \\
\mc{M}^{-1}_{10}(k,l) & \mc{M}^{-1}_{11}(k,l) \end{pmatrix} \qquad \text{for $1 \leq k,l \leq n$}.
$$
Then we have
\begin{align}
\mc{M}^{-1}_{00}(k,l) 
=&   \prod_{j=1}^n \( \frac{1- u_k v_j}{1+u_k v_j} \frac{1-u_l v_j}{1+ u_l v_j} \) 
\prod_{\begin{subarray}{c} 1 \leq i \leq n, \\ i \not= k \end{subarray}} 
\( \frac{u_k+u_i}{u_k-u_i} \)
 \prod_{\begin{subarray}{c} 1 \leq j \leq n, \\ j \not= l \end{subarray}} \( \frac{u_l + u_j}{u_l-u_j} \)
 \frac{u_k-u_l}{u_k+u_l}; \label{MI00} \\
\mc{M}^{-1}_{01}(k,l) 
=& - \mc{M}^{-1}_{10} (l,k) \label{MI01} \\ 
=&   \prod_{j=1}^n \( \frac{1- u_k v_j}{1+u_k v_j} \frac{1-u_j v_l}{1+ u_j v_l} \) 
\prod_{\begin{subarray}{c} 1 \leq i \leq n, \\ i \not= k \end{subarray}} 
\( \frac{u_k+u_i}{u_k-u_i} \)
\prod_{\begin{subarray}{c} 1 \leq j \leq n, \\ j \not= l \end{subarray}}\( \frac{v_l + v_j}{v_l-v_j} \)
 \frac{1+u_k v_l}{1-u_k v_l}; \notag \\
\mc{M}^{-1}_{11}(k,l) 
=&   -\prod_{j=1}^n \( \frac{1- u_j v_k}{1+u_j v_k} \frac{1-u_j v_l}{1+ u_j v_l} \) 
\prod_{\begin{subarray}{c} 1 \leq i \leq n, \\ i \not= k \end{subarray}} 
\( \frac{v_k+v_i}{v_k-v_i} \)
\prod_{\begin{subarray}{c} 1 \leq j \leq n, \\ j \not= l \end{subarray}} \( \frac{v_l + v_j}{v_l-v_j} \)
 \frac{v_k-v_l}{v_k+v_l}. \label{MI11}
\end{align}
\end{prop}

We recall some formulas for pfaffians, which are obtained in \cite{IW1}.

\begin{lem}[\cite{IW1}] \label{Lem-SUM}
Let $n$ be any even number.
Let $A$ and $B$ be $n$ by $n$ skew-symmetric matrices.
Then 
$$
\pf(A+B)= \sum_{\begin{subarray}{c} \bm{i} \subset \{1, \dots, n \}, \\ \text{$\# \bm{i}$: even} \end{subarray}}
(-1)^{|\bm{i}|+\frac{\# \bm{i}}{2} } \pf(A_{\bm{i}}) \pf(B_{\bm{i}^c}),
$$
where we put
$|\bm{i}|=i_1+ i_2+ \cdots + i_p$ and $\bm{i}^c = \{1, \dots, \hat{i_1}, \dots, \hat{i_2}, \dots, \hat{i_p}, \dots,
n \}$
for an increasing sequence $\bm{i}=(i_1 < i_2 < \dots < i_p)$. 
\qed
\end{lem} 

\begin{lem}[Minor summation formula \cite{IW1}] \label{MSF}
Let $m \leq n$.
Let $B$ be any $n$ by $n$ skew-symmetric matrix and $T$ be any $m$ by $n$ matrix.
Then 
$$
\sum_{\begin{subarray}{c} \bm{i} \subset \{1, \dots, n \}, \\ \# \bm{i}=m \end{subarray}}
\pf(B_{\bm{i}}) \det (T^{1,\dots,m }_{\bm{i}})
=\pf(T B {^t T}),
$$
where $T^{1,\dots,m }_{\bm{i}} = (t_{k, i_j})_{1 \leq k \leq m, 1 \leq j \leq m}$ 
for $T=(t_{ij})_{1 \leq i \leq m, 1 \leq j \leq n}$ and $\bm{i}=(i_1< \dots < i_m)$.\qed
\end{lem}

\begin{lem} \label{PfV}
Let $A$ be any $2n$ by $2n$ skew-symmetric matrix and $B$ be any $2n$ by $2m$ matrix.
Assume that $n$ is even. 
Then we have
$$
\pf \begin{pmatrix}  A & B \\ -^t B & J_{m} \end{pmatrix} = \pf(B J_m {^t B}-A).
$$
\end{lem}

\begin{proof}
By Lemma \ref{Lem-SUM} we have
\begin{align*}
\mr{LHS} =& \pf 
\( \begin{pmatrix}  A & O \\ O & J_{m} \end{pmatrix} +\begin{pmatrix}  O & B \\ -^t B & O \end{pmatrix} \) \notag \\
=& \sum_{\begin{subarray}{c} \bm{i} \subset \{1, \dots ,2n\}, \\ \text{$\# \bm{i}$: even} \end{subarray}}
\sum_{\begin{subarray}{c} \bm{k} \subset \{1, \dots ,2m\}, \\ \text{$\# \bm{k}$: even} \end{subarray}} 
(-1)^{|\bm{i}| + |\bm{k}| + (\# \bm{i} + \# \bm{k})/2} 
\pf(A_{\bm{i}}) \pf ((J_m)_{\bm{k}}) 
\pf \( \begin{pmatrix}  O & B \\ -^t B & O \end{pmatrix}_{(\bm{i}^c, n+\bm{k}^c)} \) \notag  \\
=& \sum_{\begin{subarray}{c} \bm{i} \subset \{1, \dots ,2n\}, \\ \text{$\# \bm{i}$: even} \end{subarray}}
\sum_{\begin{subarray}{c} \bm{k} \subset \{1, \dots ,2m\}, \\ \#\bm{k}^c =\#\bm{i}^c \end{subarray}} 
(-1)^{|\bm{i}| + |\bm{k}|+(\# \bm{i} + \# \bm{k})/2 + \# \bm{i}^c/2} 
\pf(A_{\bm{i}}) \pf ((J_m)_{\bm{k}}) \det(B^{\bm{i}^c}_{\bm{k}^c}).
\end{align*}
Since it is easy to see that $(\# \bm{i}+ \#\bm{i}^c)/2 = n \equiv 0 \pmod 2$ and
that $\pf ((J_m)_{\bm{k}}) = 0$ unless $\bm{k}=(2r_1-1, 2r_1, \dots, 2r_p-1, 2r_p)$ 
for an increasing sequence $(r_1, \dots, r_p)$,
we have
\begin{equation*}
\mr{LHS} = \sum_{\begin{subarray}{c} \bm{i} \subset \{1, \dots ,2n\}, \\ \text{$\# \bm{i}$: even} \end{subarray}}
\sum_{\begin{subarray}{c} \bm{k} \subset \{1, \dots ,2m\}, \\ \#\bm{k}^c =\#\bm{i}^c \end{subarray}} 
(-1)^{|\bm{i}|} 
\pf(A_{\bm{i}}) \pf ((J_m)_{\bm{k}}) \det(B^{\bm{i}^c}_{\bm{k}^c}).
\end{equation*}

On the other hand, by Lemma \ref{Lem-SUM} and Lemma \ref{MSF},
\begin{align*}
\mr{RHS} 
=& \sum_{\begin{subarray}{c} \bm{i} \subset \{1, \dots ,2n\}, \\ \text{$\# \bm{i}$: even} \end{subarray}}
(-1)^{|\bm{i}| + \frac{\# \bm{i}}{2} } \pf (-A_{\bm{i}}) \pf ((B J_m {^t B})_{\bm{i}^c}) \notag \\
=& \sum_{\begin{subarray}{c} \bm{i} \subset \{1, \dots ,2n\}, \\ \text{$\# \bm{i}$: even} \end{subarray}}
(-1)^{|\bm{i}| }\pf (A_{\bm{i}}) 
\sum_{\begin{subarray}{c} \bm{j} \subset \{1, \dots ,2m\}, \\ \#\bm{j} =\#\bm{i}^c \end{subarray}} 
\pf ((J_m)_{\bm{j}}) \det ( B^{\bm{i}^c}_{\bm{j}}).
\end{align*}
Hence we obtain the lemma.
\end{proof}

\begin{proof}[Proof of Proposition \ref{Prop-MI}]
It follows from Lemma \ref{PfV} and Proposition \ref{Prop-L} that 
\begin{align}
\pf (\mc{M}) =& \pf (J[\mf{Y}]+ L) 
= \sum_{X \in \mf{P}^{\mr{even}}(\mf{Y})} \pf (L[X \cup \{1, \dots, n \}]) \notag \\
=& \sum_{\lambda \in \mc{D}} 2^{-\ell(\lambda)} Q_{\lambda}(u) Q_{\lambda}(v) 
\prod_{1 \leq i < j \leq n} \(\frac{u_i-u_j}{u_i+u_j} \)
\prod_{1 \leq i < j \leq n} \(\frac{v_i-v_j}{v_i+v_j} \)  \notag \\
=&  \prod_{1 \leq i,j \leq n} \(\frac{1+ u_i v_j}{1-u_i v_j} \)
\prod_{1 \leq i < j \leq n} \(\frac{u_i-u_j}{u_i+u_j} \)
\prod_{1 \leq i < j \leq n} \( \frac{v_i-v_j}{v_i+v_j} \). \label{PfM}
\end{align}
We write 
$$
\mc{M}(k,l) 
= \begin{pmatrix} \mc{M}_{00} (k,l) & \mc{M}_{01} (k,l) \\
\mc{M}_{10} (k,l) & \mc{M}_{11} (k,l) \end{pmatrix}
= \begin{pmatrix}  \frac{u_k -u_l}{u_k+u_l} & - \frac{1+u_k v_l}{1-u_k v_l} \\ 
\frac{1+v_k u_l}{1- v_k u_l} & -\frac{v_k-v_l}{v_k+v_l} 
\end{pmatrix}  
$$
and $\mc{M}_{ij} = (\mc{M}_{ij}(k,l))_{1 \leq k,l \leq n}$ for $i,j=0,1$.
The skew-symmetric matrix 
$$
\widetilde{\mc{M}} = \begin{pmatrix} \mc{M}_{00} & \mc{M}_{01} \\ \mc{M}_{10} & \mc{M}_{11}
\end{pmatrix}
$$
is obtained by changing orders of rows and columns from $\mc{M}$ and it satisfies
$\pf (\widetilde{\mc{M}}) = (-1)^{\frac{n}{2}} \pf( \mc{M})$.
The pfaffian of the matrix $\widetilde{\mc{M}}^{[k,l]}$
obtained by removing the $k$-th row, $l$-th row, $k$-th column and $l$-th column
from $\widetilde{\mc{M}}$,
where $1 \leq k,l \leq n$, equals 
the value removed variables $u_k$ and $u_l$ from \eqref{PfM}, up to sign.
Namely, we have
$$
\pf (\widetilde{\mc{M}}^{[k,l]}) = (-1)^{\frac{n}{2}} 
 \prod_{\begin{subarray}{c} 1 \leq i,j \leq n, \\ i \not=k,l \end{subarray}}\( \frac{1+ u_i v_j}{1-u_i v_j} \)
\prod_{\begin{subarray}{c} 1 \leq i < j \leq n, \\ i,j \not=k,l \end{subarray}} \(\frac{u_i-u_j}{u_i+u_j} \)
\prod_{1 \leq i < j \leq n} \(\frac{v_i-v_j}{v_i+v_j} \). 
$$
Thus for $k <l$
\begin{align*}
& \mc{M}^{-1}_{00}(k,l) = (\widetilde{\mc{M}}^{-1})_{k,l}=
\frac{(-1)^{k+l} \pf (\widetilde{\mc{M}}^{[k,l]})}{ \pf (\widetilde{\mc{M}})} \\
=&  \prod_{j=1}^n \(\frac{1-u_k v_j}{1+u_k v_j} \frac{1-u_l v_j}{1+u_l v_j} \)
 \prod_{ \begin{subarray}{c} 1 \leq i \leq n, \\ i \not=k \end{subarray}}
\(\frac{u_k+u_i}{u_k-u_i} \)
 \prod_{ \begin{subarray}{c} 1 \leq j \leq n, \\ j \not=l \end{subarray}}
\(\frac{u_l+u_j}{u_l-u_j} \) \frac{u_k-u_l}{u_k+u_l}.
\end{align*}
This equality holds also if $k \geq l$ so that we obtain \eqref{MI00}.

Similarly, since $\widetilde{\mc{M}}^{[k,n+l]}$ is the matrix removed variables $u_k$ and $v_l$ 
from $\widetilde{\mc{M}}$ we see that
$$
\pf (\widetilde{\mc{M}}^{[k,n+l]}) = (-1)^{\frac{n}{2}} 
\prod_{\begin{subarray}{c} 1 \leq i,j \leq n, \\ i \not=k, j \not=l \end{subarray}} \( \frac{1+ u_i v_j}{1-u_i v_j} \)
\prod_{\begin{subarray}{c} 1 \leq i < j \leq n, \\ i,j \not=k \end{subarray}} \(\frac{u_i-u_j}{u_i+u_j} \)
\prod_{\begin{subarray}{c} 1 \leq i < j \leq n, \\ i,j \not=l \end{subarray}}\( \frac{v_i-v_j}{v_i+v_j} \). 
$$
Thus for $1 \leq k,l \leq n$ we have
\begin{align*}
& \mc{M}^{-1}_{01}(k,l) = -\mc{M}^{-1}_{10}(l,k)= (\widetilde{\mc{M}}^{-1})_{k,n+l}=
\frac{(-1)^{k+l} \pf (\widetilde{\mc{M}}^{[k,n+l]})}{ \pf (\widetilde{\mc{M}})} \\
=&   \prod_{j=1}^n \(\frac{1-u_k v_j}{1+u_k v_j} \frac{1-u_j v_l}{1+u_j v_l} \)
 \prod_{ \begin{subarray}{c} 1 \leq i \leq n, \\ i \not=k \end{subarray}}
\(\frac{u_k+u_i}{u_k-u_i} \)
 \prod_{ \begin{subarray}{c} 1 \leq j \leq n, \\ j \not=l \end{subarray}}
\(\frac{v_l+v_j}{v_l-v_j} \) \frac{1+u_k v_l}{1-u_k v_l},
\end{align*}
and so we obtain \eqref{MI01}.
We also have \eqref{MI11} in a similar way.
\end{proof}

%
\subsection{Proof of Theorem \ref{THM-SS}: final step} 
%

We calculate $K=-\eta^{-\frac{1}{2}} {^t \mc{W}} \mc{M}^{-1} \mc{W} \eta^{-\frac{1}{2}}$.
For $r,s \geq 1$ we have
\begin{equation} \label{KernelK}
K(r,s) = \begin{pmatrix} K_{00}(r,s) & K_{01}(r,s) \\ K_{10}(r,s) & K_{11}(r,s) \end{pmatrix}
=2 \sum_{k,l=1}^n \begin{pmatrix} -u_k^r \mc{M}^{-1}_{00}(k,l)u_l^s & u_k^r \mc{M}^{-1}_{01}(k,l)v_l^s \\
v_k^r \mc{M}^{-1}_{10}(k,l)u_l^s & -v_k^r \mc{M}^{-1}_{11}(k,l)v_l^s \end{pmatrix}.
\end{equation}
We prove that $\mc{K}(r,s) = K(r,s)$ for $r,s \geq 1$, where $\mc{K}$ is defined in Theorem \ref{THM-SS}.
Assume that variables $u_1, \dots, u_n, v_1, \dots, v_n$ are $2n$
distinct complex numbers in the unit open disc.

By the residue theorem
we express $\mc{K}_{00}(r,s)$ as
$$
\mc{K}_{00}(r,s) = \frac{1}{2} \frac{1}{(2\pi \sqrt{-1})^2} 
\iint_{|z|>|w|>1} \frac{Q_u(z) Q_u(w)}{Q_v(z^{-1})Q_v(w^{-1})}
\frac{z-w}{z+w} \frac{dz dw}{z^{r+1} w^{s+1}},
$$
where the contours in the integral are two circles around $0$ near the unit circle
such that $|z|>|w|>1$.
But since the integral for the residue at $z=-w$ equals
$$
\frac{(-1)^{r}}{2\pi \sqrt{-1}} \int  \frac{d w}{w^{r+s+1}} =0,
$$
we may let the contours $|z|=|w|=1$.
Changing variables $z \to z^{-1}$ and $w \to w^{-1}$,
we have
$$
\mc{K}_{00}(r,s) = \frac{1}{2} \frac{1}{(2\pi \sqrt{-1})^2} \iint_{|z|=|w|=1}
\prod_{i=1}^n \( \frac{z+u_i}{z-u_i} \frac{1-v_i z}{1+v_i z} \frac{w+u_i}{w-u_i} \frac{1-v_i w}{1+v_i w} \)
\( - \frac{z-w}{z+w} \) z^{r-1} w^{s-1} dz dw.
$$
The integrand has simple poles at $z= u_k$ and $w=u_l$ for $k,l=1, \dots, n$.
Thus it follows from \eqref{KernelK} and Proposition \ref{Prop-MI} that
\begin{align*}
\mc{K}_{00}(r,s) =& -2 \sum_{k,l=1}^n u_k^r u_l^s
 \prod_{j=1}^n \(\frac{1-u_k v_j}{1+u_k v_j} \frac{1-u_l v_j}{1+u_l v_j} \)
 \prod_{ \begin{subarray}{c} 1 \leq i \leq n, \\ i \not=k \end{subarray}}
\(\frac{u_k+u_i}{u_k-u_i} \)
 \prod_{ \begin{subarray}{c} 1 \leq j \leq n, \\ j \not=l \end{subarray}}
\(\frac{u_l+u_j}{u_l-u_j} \) \frac{u_k-u_l}{u_k+u_l} \\
=& K_{00}(r,s).
\end{align*}

Similarly, we have
$$
\mc{K}_{01}(r,s) 
= \frac{1}{2} \frac1{(2\pi \sqrt{-1})^2} \iint \frac{Q_u(z) Q_v(w)}{Q_v(z^{-1}) Q_u(w^{-1})} \frac{zw+1}{zw-1} 
\frac{d z d w}{z^{r+1} w^{s+1}}.
$$
Here the contours are two circles such that $|z|<1,\ |w|<1$.
Thus as like the calculation for $\mc{K}_{00}$ we have
\begin{align*}
\mc{K}_{01}(r,s) 
=& \frac{1}{2} \frac1{(2\pi \sqrt{-1})^2} \iint 
\prod_{i=1}^n \( \frac{z+u_i}{z-u_i} \frac{1-v_i z}{1+v_i z} \frac{w+v_i}{w-v_i} \frac{1-u_i w}{1+u_i w} \)
\frac{1+zw}{1-zw} z^{r-1} w^{s-1} dz dw \\
=& 2 \sum_{k,l=1}^n u_k^r v_l^s
 \prod_{j=1}^n \(\frac{1-u_k v_j}{1+u_k v_j} \frac{1-u_j v_l}{1+u_j v_l} \)
 \prod_{ \begin{subarray}{c} 1 \leq i \leq n, \\ i \not=k \end{subarray}}
\(\frac{u_k+u_i}{u_k-u_i} \)
 \prod_{ \begin{subarray}{c} 1 \leq j \leq n, \\ j \not=l \end{subarray}}
\(\frac{v_l+v_j}{v_l-v_j} \) \frac{1+u_k v_l}{1-u_k v_l} \\
=& K_{01}(r,s). 
\end{align*}

We can prove $\mc{K}_{11}(r,s) = K_{11}(r,s)$ in a similar way.
Hence we obtain $\mc{K}(r,s)=K(r,s)$.
Though
we have assumed that $u_i$, $v_j$ belong to the unit open disc,
it is in fact unnecessary,
see e.g. \cite{BR}.
This completes the proof of Theorem \ref{THM-SS}.

%
\section{Positivity of pfaffians} \label{Positivity}
%

We study the sufficient condition that the pfaffian is non-negative.

Let $\bSA_{2m}$ be the set of all $2m$ by $2m$ skew-symmetric $\bC$-matrices $B$ satisfying 
$\cB_{00}(r,s)=\overline{\cB_{11}(r,s)}$ and $\cB_{01}(r,s)=-\overline{\cB_{10}(r,s)}$
for each block $\cB(r,s) \ (1 \leq r,s \leq m)$.
Equivalently, $\bSA_{2m}$ 
is the set of all $2m$ by $2m$ skew-symmetric $\bC$-matrices $B$ such that $-J_m B$ is hermitian,
i.e., $B J_m =J_m \overline{B}$.
Put
$$
\bSA_{2m}^{\geq 0} = \{ B \in \bSA_{2m} ;\  \text{$-J_mB$ is non-negative definite} \}.
$$
Similarly, let $\bSA_{2m}^{> 0}$ be all elements in $\bSA_{2m}$ such that $-J_m B$ is positive definite.

\begin{prop} \label{PfaffianPositive}
Let $B \in \bSA_{2m}^{\geq 0}$.
Then $\pf(B[S])$ is non-negative for any subset $S \subset \{1, \dots, m\}$.
In particular, $\pf(B) \geq 0$.
\end{prop}

\begin{proof}
Since submatrix $-J_n B[S]$ ($S \subset \{1, \dots, m\}, \ \# S=n$) 
of $-J_m B$ is a non-negative definite hermitian matrix
by the assumption,
the roots of the polynomial $\det(zI_{2n}+J_n B[S])$ are all non-negative real numbers.
Since $\pf(z J_n -B[S])^2= \det(zI_{2n} +J_n B[S])$,
the roots of $\pf(z J_n -B[S])$ are also non-negative.
Hence if we look as the constant term in the polynomial above
we obtain that $\pf(B[S])$ is non-negative.
\end{proof}

Now we state the relation between elements in $\bSA_{2m}^{\geq 0}$ and non-negative definite hermitian matrices. 
Let $\mbf{Her}_m$ be the set of all $m$ by $m$ hermitian matrices.
Let $\mbf{Her}_m^{\geq 0} \subset \mbf{Her}_m$ be  all non-negative definite hermitian matrices.

\begin{prop}
Let $\omega:\fgl_m(\bC) \to \fo_{2m}(\bC)$ be an injective linear map defined in Section \ref{Connection}.
Then we have $\omega(\mbf{Her}_m) \subset \bSA_{2m}$ and $\omega(\mbf{Her}_m^{\geq 0}) \subset \bSA_{2m}^{\geq 0}$.
\end{prop}

\begin{proof}
For $A=(a_{rs})_{1 \leq r,s \leq m} \in \mbf{Her}_m$ the $(r,s)$-block of $\omega(A)$ is 
$$
\omega(A)(r,s)= \begin{pmatrix} 0 & a_{rs} \\ -\overline{a_{rs}} & 0 \end{pmatrix}
$$
so that $\omega(A) \in \bSA_{2m}$.
Let $\lambda$ be an eigenvalue of $-J_m \omega(A)$.
Then since 
\begin{align*}
0=& \det (\lambda I_{2m} +J_m \omega(A))=
\pf (\lambda J_m -\omega(A))^2 
= \det \begin{pmatrix} O_m & \lambda I_m-A \\ -\lambda I_m +\Bar{A} & O_m \end{pmatrix} \\
= &\det(\lambda I_m-A) \det(\lambda I_m-A^*)
\end{align*}
and $A$ and $A^*$ are non-negative definite,
the eigenvalue $\lambda$ is non-negative.
\end{proof}

Hence we may treat elements in
$\bSA_{2m}$ as generalizations of hermitian matrices.
Proposition \ref{PfaffianPositive} indicates us 
the following problem.

\medskip

\begin{problem} \label{pro}
What is the range for $\alpha$ such that $\pfa{\alpha}(B) \geq 0$
whenever $B \in \bSA_{2m}^{\geq 0}$.
\qed
\end{problem}

\medskip

This problem is an extension of the problem for the positivity of 
$\alpha$-determinants of a non-negative definite hermitian matrix
(recall Conjecture \ref{ConjST}).
It relates to the existence of $\alpha$-pfaffian point processes defined in Section \ref{PPP}.

\begin{prop}
Let $L$ be a skew-symmetric matrix kernel on a countable set $\mf{X}$ and
$\alpha$ be a real number.
Assume that $-JL$ is hermitian and non-negative definite and satisfies $\norm{\alpha L}<1$.
Then $-JK_\alpha$ is also hermitian and non-negative definite,
where $K_\alpha=L(I+\alpha JL)^{-1}=(I+\alpha LJ)^{-1} L$.
\end{prop}

\begin{proof}
Put $\mc{L}=-JL$ and $\mc{K}_\alpha= -JK_\alpha$.
It is easy to see that $\mc{K}_\alpha= \mc{L}(I-\alpha \mc{L})^{-1}= (I-\alpha \mc{L})^{-1} \mc{L}$
so that $\mc{K}_\alpha$ is hermitian if $\alpha$ is a real number.
If $\lambda$ is an eigenvalue of $\mc{L}$,
then $|\alpha \lambda|<1$ since $\norm{\alpha \mc{L}} <1$.
Therefore $\mc{K}_\alpha$ is non-negative definite.
\end{proof}

For an $\alpha$, we suppose that  
$\pfa{\alpha}(L[X])$ is non-negative for all $X \in \mf{P}(\mf{X})$
whenever $-JL$ is hermitian and non-negative definite.
Then, by the proposition above, the $\alpha$-pfaffian $\pfa{\alpha}(K_\alpha[X])$
is also non-negative.
Therefore the correlation function in Theorem \ref{THM-Correlation} exists in this case.

\medskip

We study Problem \ref{pro}.
It follows from Proposition \ref{PfaffianPositive} that
$\pfa{-1}(B)=\pf (B) \geq 0$ and $\pfa{0}(B) \geq 0$
for $B \in \bSA_{2m}^{ \geq 0}$.
More strongly, we obtain the following inequality.

\begin{thm} \label{pfa0>pf}
Let $B \in \bSA_{2m}^{>0}$.
Then we have
$$
\pfa{0}(B)= \prod_{r=1}^m \cB_{01}(r,r) \geq \pf (B).
$$
The equality holds
if and only if $\cB(r,s)=O_{2}$ for $r \not=s$.
\end{thm}

If we let $B=\omega(A)$ for $A=(a_{ij}) \in {\bf Her}_m$, then Theorem \ref{pfa0>pf} provides the well-known formula
$$
\prod_{i=1}^m a_{ii} \geq \det(A).
$$
We need some preparations to prove the theorem.
Put 
$$
\cS_{2m,2} = \left\{ X= 
\begin{pmatrix} X^{(1)} \\ \vdots \\ X^{(m)} \end{pmatrix}
;\  X^{(j)} \in \bSA_{2} \ (1 \leq j \leq m ) \right\}.
$$
It is easy to see that $J_m X =\overline{X} J_1$ for $X \in \cS_{2m,2}$.
Let $O_{r,s}$ be the $r$ by $s$ zero matrix.

\begin{lem} \label{Lem-Positivity}
For $B \in \bSA_{2m}^{>0}$ and $X \in \cS_{2m,2}\setminus \{O_{2m,2}\}$, we have $\pf( {^t X} B X )>0$.
\end{lem}

\begin{proof}
Since the matrix ${^t X}B X$ is a $2$ by $2$ skew-symmetric matrix,
it is easy to see that 
$\pf( {^t X} BX)= \frac{1}{2} \tr(-J_1 {^t X} BX)$. 
Therefore it is sufficient to prove $\tr (-J_1 {^t X} BX)>0$.
By the commutation relation $J_m X= \overline{X}J_1$, we have 
$\tr (-J_1 {^t X} BX)= \tr ( {^t \overline{X}} (-J_m B)X)$.
If we write $X=( \bm{x}_1, \bm{x}_2)$ ($\bm{x}_1, \bm{x}_2 \in \bC^{2m}$),
we have 
$\tr ( {^t \overline{X}} (-J_m B)X)
= {^t \overline{\bm{x}_1} } (-J_m B) \bm{x}_1 + {^t \overline{\bm{x}_2}} (-J_m B) \bm{x}_2 > 0$
since $-J_m B$ is positive definite and $\bm{x}_1 \not=\bm{0}$ or $\bm{x}_2 \not= \bm{0}$.
\end{proof}

\begin{lem} \label{Lem-Inverse}
If $B \in \bSA_{2m}^{>0}$, then $(B^*)^{-1} \in \bSA_{2m}^{>0}$.
\end{lem}

\begin{proof}
Let $B \in \bSA_{2m}^{>0}$.
Since $-J_m B$ is positive definite and hermitian,
the inverse $(-J_mB)^{-1}= B^{-1} J_m=J_m \overline{B}^{-1} = -J_m(B^{*})^{-1}$ is also.
Therefore $(B^*)^{-1} \in \bSA_{2m}^{>0}$.
\end{proof}

\begin{proof}[Proof of Theorem \ref{pfa0>pf}]
Write $B$ in the form
$$
B= \begin{pmatrix} B_{m-1} & X \\ -^tX & \cB(m,m) \end{pmatrix}, \qquad B_{m-1} \in \bSA_{2m-2},\ 
X = \begin{pmatrix} \cB(1,m) \\ \vdots \\ \cB(m-1,m) \end{pmatrix} \in \cS_{2m-2,2}.
$$
Then $B_{m-1} \in \bSA_{2m-2}^{>0}$. 
In fact, 
if $B \in \bSA_{2m}^{>0}$, then $-J_m B$ is positive definite.
We immediately see that $-J_{m-1} B_{m-1}$ is also positive definite, hence $B_{m-1} \in \bSA_{2m-2}^{>0}$.

Put 
$$
P= \begin{pmatrix} I_{2m-2} & (B_{m-1})^{-1} X \\ O_{2,2m-2} & I_2 \end{pmatrix}.
$$
Then it is easy to see that
$$
B= {^t P} \begin{pmatrix} B_{m-1} & O_{2m-2,2} \\ O_{2,2m-2} & \cB(m,m) + {^t X}(B_{m-1})^{-1} X \end{pmatrix} P.
$$
We have
\begin{align*}
\pf(B)=& \pf(B_{m-1}) \pf(\cB(m,m) + {^t X}(B_{m-1})^{-1} X ) \det(P) \\
=&
\pf(B_{m-1}) \{ \pf(\cB(m,m)) - \pf({^t X}(-B_{m-1})^{-1} X ) \}.
\end{align*}
Since $B_{m-1} \in \bSA_{2m-2}^{>0}$,
it follows from Lemma \ref{Lem-Positivity} and Lemma \ref{Lem-Inverse}
that
$\pf({^t \overline{X}}(B_{m-1}^*)^{-1}) \overline{X} )>0$.
Hence
$\pf ({^t X} (-B_{m-1})^{-1} X) >0$ for any $X \in \mc{S}_{2m-2,2}$
 if and only if $X \not= O_{2m-2,2}$.
Therefore we have
$$
\pf(B) \leq \pf(B_{m-1}) \pf (\cB(m,m))
$$
and the equality holds if and only if $X=O_{2m-2,2}$.
By the induction on $m$, we have the theorem.
\end{proof}

%
%

\noindent
\textsc{MATSUMOTO, Sho}\\
Graduate School of Mathematics, Kyushu University.\\
Hakozaki Higashi-ku, Fukuoka, 812-8581 JAPAN.\\
\texttt{ma203029@math.kyushu-u.ac.jp}\\

\end{document}